%% file: weak-palis-revised.tex
\documentclass[11pt]{article}
\usepackage{amsmath,amsfonts,amsthm,epsfig,latexsym, amssymb}

 \def\NN{{\mathbb N}}  
 \def\RR{{\mathbb R}}  
   
 \def\ZZ{{\mathbb Z}}

  \def\cG{{\cal G}} \def\cM{{\cal M}} \def\cS{{\cal S}}
  \def\cH{{\cal H}}  \def\cT{{\cal T}}
\def\cC{{\cal C}}  \def\cI{{\cal I}} \def\cO{{\cal O}} \def\cU{{\cal U}}
\def\cD{{\cal D}}    
    
   \def\cR{{\cal R}}

\def\La{\Lambda}

\def\Om{\Omega}

\newcommand{\dd}{\operatorname{d}}

\newcommand{\closure}{\operatorname{Cl}}
\newcommand{\interior}{\operatorname{Int}}
\newcommand{\boundary}{\partial}

\newcommand{\vol}{{v}}

\newcommand{\diff}{\operatorname{Diff}}

\newtheorem*{maintheorem}{Main theorem}
\newtheorem*{theoremA}{Theorem A}
\newtheorem*{theoremB}{Theorem B}
\newtheorem*{theorem*}{Theorem}
\newtheorem*{weak}{Weak Palis conjecture}
\newtheorem*{strong}{Strong Palis conjecture}
\newtheorem*{plaque}{Plaque family theorem}
\newtheorem*{corollary*}{Corollary}

\newtheorem{lemma}{Lemma}[section]
\newtheorem{theorem}[lemma]{Theorem}
\newtheorem{corollary}[lemma]{Corollary}
\newtheorem{proposition}[lemma]{Proposition}

\newtheorem{addendum}[lemma]{Addendum}
\theoremstyle{definition}
\newtheorem{definition}[lemma]{Definition}
\theoremstyle{remark}
\newtheorem{remark}[lemma]{Remark}

\numberwithin{equation}{section}

\makeatletter
\renewcommand{\l@section}{\@dottedtocline{2}{3.8em}{3.2em}}
\renewcommand{\l@subsection}{\@dottedtocline{3}{3.8em}{3.2em}}
\newcommand{\subsectionruninhead}{\@startsection{subsection}{2}{0mm}{-\baselineskip}{-0mm}{\bf\large}}
\newcommand{\subsubsectionruninhead}{\@startsection{subsubsection}{3}{0mm}{-\baselineskip}{-0mm}{\bf\normalsize}}
\makeatother

\newcommand{\footnoteagain}[1]{\mbox{$^\text{#1}$}}

\begin{document}
\title{Birth of homoclinic intersections: a model for the central dynamics of partially hyperbolic systems}
\author{Sylvain Crovisier}
\date{4-12-2009}

\maketitle
\begin{abstract}
We prove a conjecture of J. Palis: any diffeomorphism of a compact manifold
can be $C^1$-approximated by a Morse-Smale diffeomorphism or by
a diffeomorphism having a transverse homoclinic intersection.

\begin{center}
\bf R\'esum\'e
\end{center}
\emph{Cr\'eation d'intersection homoclines~: un mod\`ele pour la dynamique centrale
des syst\`emes partiellement hyperboliques.}
Nous montrons une conjecture de J. Palis~: tout diff\'eomorphisme d'une vari\'et\'e
compacte peut \^etre approch\'e en topologie $C^1$ par un diff\'eomorphisme Morse-Smale
ou par un diff\'eomorphisme ayant une intersection homocline transverse.

\vskip 2mm

\begin{description}
\item[\bf Key words:] Homoclinic orbit, Morse-Smale diffeomorphism, hyperbolic diffeomorphism, generic dynamics,
homoclinic class, partial hyperbolicity.
\item[\bf MSC 2000:] 37B40, 37C05, 37C20, 37C29, 37C50, 37D15, 37D30.
\end{description}
\end{abstract}

\setcounter{section}{-1}
\tableofcontents

%%%%%%%%%%%%%%%%%%%%%%%%%%%%%%%%%%%%%%%%%%%%%%%%%%%%%%%%%%%%%%%%%%%%%%%%%%%%%%%%%%
%%%%%%%%%%%%%%%%%%%%%%%%%%%%%%%%%%%%%%%%%%%%%%%%%%%%%%%%%%%%%%%%%%%%%%%%%%%%%%%%%%
\section{Introduction}
%%%%%%%%%%%%%%%%%%%%%%%%%%%%%%%%%%%%%%%%%%%%%%%%%%%%%%%%%%%%%%%%%%%%%%%%%%%%%%%%%%
\subsection{Homoclinic intersections in dynamical systems}
In his (corrected) prize essay~\cite{poincare1}, H. Poincar\'e was the first to imagine around 1890 the existence of ``doubly asymptotic solutions",
that are now called \emph{transverse homoclinic intersections}: for a diffeomorphism,
these are the transverse intersection points between the stable and unstable manifolds of a hyperbolic periodic orbit.
He discovered that this simple assumption implies a very intricate behavior: in~\cite{poincare2}, he showed that any transverse homoclinic intersection
is accumulated by other homoclinic points. In 1935 G.D. Birkhoff proved~\cite{birkhoff} by topological arguments that it is also the limit of infinitely
many periodic points. Later on (in 1965), S. Smale introduced \cite{smale-horseshoe} a geometrical model
for the dynamics occurring near a homoclinic orbit: the horseshoe is an invariant compact set that carries
a hyperbolic dynamics and contains all the complicated phenomena discovered
by Poincar\'e and Birkhoff, but that can also be described by coding. Hence, the transverse homoclinic intersections give birth to very rich dynamics
(in particular, the topological entropy is non zero) and moreover they are robust (i.e. they persist under $C^1$-perturbation of the system).
S. Newhouse then associated \cite{newhouse-homocline} to any hyperbolic periodic orbit having transverse homoclinic intersections a very natural invariant compact set
that contains these horseshoes: he defined the \emph{homoclinic class} as the closure of all
the transverse homoclinic points of this periodic orbit.
Nowadays the phenomena creating homoclinic bifurcations still remain a major topic in differentiable dynamical systems
(see~\cite{palis-takens, bonatti-diaz-viana}).

In the opposite direction, Smale defined~\cite{smale-morse} a class of diffeomorphisms with very simple dynamics, the \emph{Morse-Smale diffeomorphisms}:
the chain recurrent set (which carries the non-trivial dynamics, see section~\ref{s.notation} for the definitions) is hyperbolic and
finite\footnote{\label{f.def-classic} The classical formulation of this
definition involves the non-wandering set (see section~\ref{s.notation}).}; moreover the invariant manifolds of the periodic orbits are in general position
(satisfying the strong transversality condition).
As a consequence, any non-periodic orbit accumulates in the past and in the future on two different periodic orbits;
these diffeomorphisms have no homoclinic intersection.
Moreover, the set of these diffeomorphisms is open:
in fact, J. Palis and Smale~\cite{palis-smale} have proven that
any $C^1$-perturbation of a Morse-Smale diffeomorphism is a Morse-Smale diffeomorphism that is conjugate to the initial dynamics through a homeomorphism
(the dynamics is structurally stable).

Palis has stated~\cite{palis-conjecture, palis-conjecture2} many conjectures about dissipative differentiable dynamics and in particular
that the homoclinic intersections are ``always" responsible for complicated dynamics:
this was formulated for the set $\diff^r(M)$ of $C^r$-diffeomorphisms ($r\geq 1$)
of a compact manifold $M$ endowed with the $C^r$-topology.

\begin{weak}
The space $\diff^r(M)$ of $C^r$-diffeomorphisms ($r\geq 1$) contains a dense open set which
decomposes as the union $\mathcal{MS}\cup\mathcal{I}$ of two disjoint open sets:
\begin{itemize}
\item $\mathcal{MS}$ is the set of Morse-Smale diffeomorphisms,
\item $\mathcal{I}$ is the set of diffeomorphisms having a transverse homoclinic intersection.
\end{itemize}
\end{weak}
One interest of the conjecture is that it gives information
on a large (open and dense) set of diffeomorphisms.

One may first restrict this question to the open set of \emph{hyperbolic diffeomorphisms}
which are the diffeomorphisms with hyperbolic chain-recurrent set\footnoteagain{\ref{f.def-classic}}.
Smale has shown~\cite{smale-horseshoe} that the chain-recurrent set of any hyperbolic diffeomorphism
is the disjoint union of finitely
many basic sets: these are the homoclinic classes defined by Newhouse in~\cite{newhouse-homocline}.
If there is no homoclinic intersection, each basic piece is trivial (it is a periodic orbit) and the dynamics is Morse-Smale (up to a $C^r$-perturbation which gives the transversality between the invariant
manifolds of the periodic orbits).
Hence, the conjecture is satisfied inside the set of hyperbolic diffeomorphisms\footnote{A more elementary argument is given at section~\ref{ss.simple}.}.

However it was discovered that the set of hyperbolic diffeomorphisms is not dense among differentiable dynamics
(two families of counter-examples were described, starting from
those built by R. Abraham and Smale~\cite{abraham-smale} and Newhouse~\cite{newhouse-thesis}).
Palis then conjectured~\cite{palis-conjecture, palis-conjecture2} that the homoclinic bifurcations are ``always" responsible for non-hyperbolicity:

\begin{strong}
Any $C^r$-diffeomorphism ($r\geq 1$) of a compact manifold can be approximated by
a diffeomorphism which is either hyperbolic, has a homoclinic tangency or has a heterodimensional cycle.
\end{strong}

It is easy to see that this second conjecture implies the first one since the dynamics at a homoclinic
bifurcation can be perturbed in order to create a transverse homoclinic intersection (see section~\ref{ss.simple}).
Note also that these conjectures hold when $M$ is one-dimensional since in this case M. Peixoto
showed~\cite{peixoto} that the Morse-Smale diffeomorphisms are dense in $\diff^r(M)$ for any $r\geq 1$.

%%%%%%%%%%%%%%%%%%%%%%%%%%%%%%%%%%%%%%%%%%%%%%%%%%%%%%%%%%
\subsection{Main results}
In the $C^r$ topologies, with $r>1$, very few techniques are available up to here for proving these conjectures
(in particular, a $C^r$-closing lemma is missing). On the other hand,
many tools have been developed during the last decade for the topology $C^1$ allowing substantial progress
in this direction.

The $C^1$ strong conjecture has been solved on surfaces
by E. Pujals and M. Sambarino~\cite{pujals-sambarino};
in higher dimensions some progress\footnote{Since the first version of this text we have obtained further results towards the strong conjecture~\cite{modele,cp}.} have been obtained by several people:
S. Hayashi, Pujals~\cite{pujals1,pujals2} and L. Wen~\cite{wen}.
Very recently, C. Bonatti, S. Gan and Wen have managed to show \cite{bonatti-gan-wen}
the weak conjecture on three-dimensional manifolds. In this paper we are aimed to prove it in
any dimension.

\begin{maintheorem}
Any diffeomorphism of a compact manifold can be $C^1$-approximated by a Morse-Smale diffeomorphism
or one exhibiting a transverse homoclinic intersection.
\end{maintheorem}

As it was noticed in~\cite{pujals-sambarino}, this result gives a
description of the set of diffeomorphisms having zero topological entropy:
we will say that a diffeomorphism $f$ has \emph{robustly
zero entropy} if any diffeomorphism in a $C^1$-neighborhood of $f$
has zero entropy. We define in a similar way the diffeomorphisms having robustly non-zero entropy.
Note that if a diffeomorphism is Morse-Smale it has zero entropy, if
it has a transverse homoclinic intersection it has non-zero entropy. As a consequence,
we get:

\begin{corollary*}
The closure of the set of diffeomorphisms having robustly zero entropy
coincides with the closure of the Morse-Smale diffeomorphisms.

The closure of the set of diffeomorphisms having robustly positive entropy
coincides with the closure of the set of diffeomorphisms exhibiting a transverse homoclinic
intersection.
\end{corollary*}

This can be compared to a property proved by A. Katok \cite{katok}
for $C^2$-diffeomorphisms on surfaces:
\emph{the $C^2$ surface diffeomorphisms having positive entropy are the diffeomorphisms exhibiting a transverse
homoclinic intersection.}

\subsection{Related results}
The $C^1$-strong conjecture is known to be true
(with a better conclusion)
among diffeomorphisms whose dynamics splits into finitely many pieces only:
these are diffeomorphisms which don't admits \emph{filtrations}\footnote{A diffeomorphism
$f$ admits \emph{filtrations with an arbitrarily large number of levels}
if there exists arbitrarily long sequences of
open sets $U_0\subset U_1\subset \dots\subset U_N\subset M$
satisfying $f(\overline{U_i})\subset U_i$
and such that the maximal set inside $U_i\setminus \overline{U_{i+1}}$
is non-empty for each $0\leq i<N$.
Equivalently $f$ has infinitely many chain-recurrence classes,
see section~\ref{s.techniques}. A diffeomorphism which can not be
$C^1$ approximated by a diffeomorphism admiting filtrations with an arbitrarily large number of levels is sometimes called \emph{tame}.
} with an arbitrarily large number of levels.

\begin{theorem*}
Any diffeomorphism of a compact manifold can be $C^1$-approximated
by a diffeomorphism which is hyperbolic or has a heterodimensional cycle
or admits filtrations with an arbitrarily large number of levels.
\end{theorem*}
The result was proved on surfaces by Ma\~n\'e~\cite{mane}.
With the results of~\cite{BC},
this theorem is just an improvement from~\cite[theorem C]{flavio} and~\cite[theorem D]{gan-wen}.
The main argument is based on Liao's and Ma\~n\'e's techniques~\cite{liao2,mane}.
This is done in the appendix.
\bigskip

Let us also mention some related results in the conservative setting. If $M$ is endowed with a
symplectic form $\omega$ or a volume form $\vol$,
one considers the spaces $\diff^r(M,\omega)$ and $\diff^r(M,\vol)$ of $C^r$ diffeomorphisms 
that preserve $\omega$ or $\vol$. Since there is no conservative Morse-Smale dynamics, one
expects that the set of conservative diffeomorphisms having a transverse homoclinic intersection is
dense.
\begin{itemize}
\item[-] For the $C^r$ topology, $r\geq 1$, some results are known on surfaces:
E. Zehnder has proven~\cite{zehnder} that any surface diffeomorphism having an elliptic periodic point is approximated by diffeomorphisms having a transverse homoclinic orbit.
More generally, the existence of transverse homoclinic intersection is dense in most homotopy class of conservative surface $C^r$-diffeomorphisms: see the works by C. Robinson~\cite{robinson-homocline},
D. Pixton~\cite{pixton} and F. Oliveira~\cite{oliveira1,oliveira2}.
\item[-] For the $C^1$-topology, one has a better
description: by the conservative closing lemma of C. Pugh and Robinson~\cite{pugh-robinson},
for any diffeomorphism $f$ in a dense G$_\delta$ subset of $\diff^1(M,\omega)$ or $\diff^1(M,\vol)$,
the hyperbolic periodic points are dense in $M$.
Moreover, F. Takens has proven~\cite{takens} (see also~\cite{xia})
that each hyperbolic periodic point
has a transverse homoclinic orbit.
This gives the $C^1$ weak Palis conjecture in the conservative setting.

\item[-] With M.-C. Arnaud and Bonatti, we have shown~\cite{BC,ABC}
that if $M$ is connected it is $C^1$-generically a single homoclinic class,
implying that the dynamics of $f$ is transitive
(there exists an orbit that is dense in $M$).
One can thus argue as for diffeomorphisms
which don't admits filtrations with an arbitrarily large number of levels. The following is proved in the appendix.

\end{itemize}

\begin{theorem*}
Any conservative diffeomorphism of a compact surface can be $C^1$-approximated
by a conservative diffeomorphism which is Anosov or has a homoclinic tangency.

Any conservative diffeomorphism of a manifold $M$ of dimension
$\dim(M)\geq 3$ can be $C^1$-approximated by a conservative diffeomorphism which is Anosov or has a heterodimensional cycle.
\end{theorem*}

%%%%%%%%%%%%%%%%%%%%%%%%%%%%%%%%%%%%%%%%%%%%%%%%%%%%%%%%%%%%%%%%%%%%%%%%%%%%%%%%%%%%%%%%%%%%%%%%%%%
\subsection{Discussions on the techniques}\label{s.techniques}
One possible approach for analyzing a dissipative dynamical system
is to \emph{split the dynamics into elementary pieces}
that one could study independently.
This method is now well understood for the $C^1$-generic diffeomorphisms:
by using pseudo-orbits, C. Conley's theory~\cite{conley} decomposes the dynamics on the chain-recurrent set
of any homeomorphism into elementary pieces, called \emph{chain-recurrence classes}.
In particular, if each chain-recurrence class of a generic diffeomorphism is reduced to a periodic orbit,
then the dynamics is Morse-Smale.
Following ideas developed by Pugh for the closing lemma~\cite{pugh} and by
Hayashi for the connecting lemma~\cite{hayashi}, we proved with Bonatti a connecting lemma for pseudo-orbits~\cite{BC}: for the $C^1$-generic diffeomorphisms it implies
that a chain-recurrence class which contains a periodic point $P$ coincides with the homoclinic class of $P$.
In order to prove the weak conjecture, it thus remains to study the case there also exist chain-recurrence classes which do not contain any periodic orbit (these are called \emph{aperiodic classes}).

We have shown~\cite{crovisier} that, for $C^1$ generic diffeomorphisms, the aperiodic classes
are the limit for the Hausdorff distance of a sequence of hyperbolic periodic orbits.
If one shows that the invariant manifolds have uniform size, one can expect that they intersect,
implying the existence of a homoclinic intersection.
For non-hyperbolic diffeomorphisms,
we will try to \emph{get weaker forms of hyperbolicity}:
we look for the existence of \emph{dominated splittings} or of \emph{partially hyperbolic structures}.
One can hope to obtain such properties if one works with diffeomorphisms
that are far from the known obstructions to the uniform hyperbolicity.

Once we are reduced to analyzing partially hyperbolic systems, the size of the invariant manifold is
controlled in the directions where the dynamics is uniformly contracting or expanding.
It remains to understand what happens in the other directions: this is much simpler if
\emph{the weak directions are one-dimensional}, since the order of $\RR$ can be used in this case.\bigskip

Let us explain how these ideas already appear in the previous works.
On surfaces, Pujals and Sambarino~\cite{pujals-sambarino} obtain a dominated splitting $E\oplus F$
for diffeomorphisms that are $C^1$-far from homoclinic tangencies.
Then, they apply $C^2$ Denjoy techniques to control the distortion,
the size of the invariant manifold and the expansion or contraction in both directions.
They strongly use the fact that each weak bundle $E$ or $F$ is one-dimensional
but also that they are extremal: it is possible to generalize their method
in higher dimension in order to analyze partially hyperbolic dynamics
with a one-dimensional central direction but under the assumption
that the strong stable or the strong unstable direction is trivial (see~\cite{pujals-sambarino-codimension} and~\cite{zhang}).

In higher dimensions, Bonatti, Gan and Wen have managed~\cite{bonatti-gan-wen}
to control a non-extremal central bundle but their argument
is semi-global and the partially hyperbolic sets they consider have to be chain-recurrent classes:
they focus on the periodic orbits and the existence of filtrating neighborhoods implies
long central stable or unstable manifolds.
In dimension three, they can prove that any $C^1$-generic diffeomorphism far from
homoclinic bifurcations has such a partially hyperbolic chain-recurrence class.

In dimension larger or equal to $4$, their argument is not sufficient since one would have to consider
chain-recurrence classes with mixed behaviors, i.e. containing partially
hyperbolic proper subsets whose stable directions have different dimension.
This prevents the existence of a partially hyperbolic structure on the whole chain-recurrence class.
\bigskip

In this paper we give a different method which allows to deal with any partially hyperbolic
chain-transitive set and to perform local arguments.
As in the previous works, the proof of the main theorem is obtained in two steps:

\begin{theoremA}
There exists a dense G$_\delta$ subset $\cG_{A}\subset \diff^1(M)$
such that for any diffeomorphism $f\in \cG_{A}\setminus (\cM\cS\cup \cI)$,
there is an invariant  compact set $K_{0}$ which carries a minimal dynamics
but is not a periodic orbit
and which is partially hyperbolic with a one-dimensional central bundle.
\end{theoremA}

\begin{theoremB}
There exists a dense G$_\delta$ subset $\cG_B\subset \diff^1(M)$ such that
any diffeomorphism $f$ in $\cG_B$ and any chain-transitive set $K_0$ of $f$
which is partially hyperbolic with a one-dimensional central bundle satisfy the following property:
If $K_0$ is not a periodic orbit, then it is the limit for the Hausdorff distance
of non-trivial relative homoclinic classes.
\end{theoremB}

These two theorems are proven in the last section of the paper. Theorem A is a consequence
of Wen's study~\cite{wen} of diffeomorphisms $C^1$-far from homoclinic bifurcations and
of the minimally non-hyperbolic sets.
Theorem B is obtained by introducing \emph{central models}: they
are defined and studied in an abstract way in section~\ref{s.abstract};
then in section~\ref{s.application} we show how they appear in differentiable dynamics.
In particular we explain how they imply the existence
of homoclinic intersections: some of the arguments here were already used
in~\cite{bonatti-gan-wen}.

A difficulty when one wants to study a partially hyperbolic system
is that we don't know in general if an invariant central foliation exist.
M. Hirsch, Pugh and M. Shub~\cite{HPS} however have shown that
a weaker property always holds: there is a locally invariant family of central plaques.
This object has been already used many times (see for instance~\cite{pujals-sambarino}
or the ``fake foliations" in the work by K. Burns and A. Wilkinson~\cite{burns-wilkinson}).
The central model is a tool for understanding the dynamics along the plaques;
in particular, we forget that different central plaques may have wild intersections.
The second ingredient is to replace the $C^2$ distortion techniques in~\cite{pujals-sambarino}
by topological arguments. In particular, one proves topological contraction in the central direction
by looking for attractors in the central models: they are obtained by using pseudo-orbits and
arguments of Conley's theory.

\vskip 5mm

\noindent
{\large \bf Acknowledgments:} {\it I am grateful to Christian Bonatti for
many comments that simplified some parts of the proof. I also want to thank
F. Abdenur, S. Gan, J. Palis, E. Pujals, M. Sambarino, L. Wen, A. Wilkinson,
for discussions during the preparation of the paper.\bigskip}

%%%%%%%%%%%%%%%%%%%%%%%%%%%%%%%%%%%%%%%%%%%%%%%%%%%%%%%%%%%%%%%%%%%
%%%%%%%%%%%%%%%%%%%%%%%%%%%%%%%%%%%%%%%%%%%%%%%%%%%%%%%%%%%%%%%%%%%
\section{Notations and definitions}\label{s.notation}
%%%%%%%%%%%%%%%%%%%%%%%%%%%%%%%%%%%%%%%%%%%%%%%%%%%%%%%%%%%%%%%%%%%

The closure, the interior and the boundary of a subset $Y$ of a topological space $X$
will be denoted by $\closure(Y)$, $\interior(Y)$ and $\boundary(Y)$ respectively.

For a compact metric space $(X,d)$, we will use the Hausdorff distance
between non-empty compact subsets $A$, $B$ of $X$, defined by:
$$\dd_H(A,B)=\max\left(\sup_{x\in A}\dd(x,B),\sup_{x\in B}\dd(x,A)\right).$$

We will say that a sequence $(z_n)$ in $X$
is a $\varepsilon$-\emph{pseudo-orbit} if $\dd(h(z_n),z_{n+1})$ is less than $\varepsilon$
for each $n$. For any points $x,y\in X$ we write $x\dashv y$ if for any $\varepsilon>0$,
there exists a $\varepsilon$-pseudo-orbit $(x=z_0,z_1,\dots,z_n=y)$ with $n\geq 1$.\\
The set $X$ is \emph{chain-transitive} if for any $x,y\in X$, we have $x\dashv y$.
An invariant compact subset $A\subset X$ is \emph{chain-transitive} if $h$ induces a chain-transitive
homeomorphism on $A$.\\
The union of the chain-transitive
subsets of $X$ is the \emph{chain-recurrent} set of $h$, denoted by $\cR(h)$.
The chain-transitive sets of $X$ that are maximal for the inclusion define a partition of $\cR(h)$
into disjoint invariant compact sets called \emph{chain-recurrent classes}.
They are the objects of Conley's theory (see \cite{conley, robinson-livre}).\\
Usually one also considers another kind of recurrence: the \emph{non-wandering set} $\Om(h)$ is the set of points that do not have any
neighborhood $U$ disjoint from all its iterates $f^k(U)$, $k\neq 0$.\bigskip

Let us consider a compact smooth manifold $M$ endowed with a Riemannian metric
and a diffeomorphism $f$ of $M$.

The diffeomorphism $f$ of $M$ is \emph{Kupka-Smale} if all its periodic orbits are hyperbolic
and if moreover for any periodic points $p$ and $q$,
the unstable manifold $W^{u}(p)$ of $p$ and the stable manifold $W^{s}(q)$ of $q$
are in general position (i.e. at any intersection point $x\in W^u(p)\cap W^s(q)$,
we have $T_xM=T_xW^{u}(p)+T_xW^{s}(q)$).
I. Kupka~\cite{kupka} and Smale~\cite{smale-kupka} have shown that for any $r\geq 1$
the set of Kupka-Smale diffeomorphisms
is a dense G$_\delta$ subset $\cG_{KS}$ of $\diff^r(M)$ endowed with the $C^r$-topology.

The diffeomorphism $f$ is \emph{Morse-Smale} if
 it is a Kupka-Smale diffeomorphism whose non-wandering set $\Omega(f)$
(or equivalently whose chain-recurrent set $\cR(f)$) is finite.

More generally, we will say that $f$ is
\emph{hyperbolic} if it is Axiom $A$ (i.e. its non-wandering set
is hyperbolic and contains a dense set of periodic points) and
moreover satisfies the no-cycle condition.
By Smale spectral theorem~\cite{smale-dynamics}, this is equivalent to requiring that
the chain-recurrent set $\cR(f)$ of $f$ is hyperbolic, which is the definition we will use in this paper.
The set of hyperbolic diffeomorphisms is an open set denoted by $\cH\text{yp}$.

A diffeomorphism has a \emph{homoclinic tangency} if there exists a hyperbolic periodic point $p$
whose stable and unstable manifold have a non-transverse intersection point.
The set of these diffeomorphisms will be denoted by $\cT\text{ang}$.

A diffeomorphism has a \emph{heterodimensional cycle} if there exist two hyperbolic periodic orbits
$\cO_1$ and $\cO_2$ whose stable manifolds have different dimensions
and such that the unstable manifold of $\cO_i$ intersects the stable manifold of $\cO_j$
for $(i,j)=(1,2)$ and $(i,j)=(2,1)$. The set of these diffeomorphisms will be denoted by $\cC\text{yc}$.\bigskip

We say that two hyperbolic periodic orbits
$\cO_1$ and $\cO_2$ are \emph{homoclinically related} if the unstable manifold of $\cO_i$ intersects
transversally the stable manifold of $\cO_j$ for $(i,j)=(1,2)$ and $(i,j)=(2,1)$.
(The intersection points are called \emph{heteroclinic points of $\cO_1$ and $\cO_2$}.)
This defines an equivalence relation on the set of hyperbolic periodic orbits.

Let $P$ be a hyperbolic periodic point. Its \emph{homoclinic class} $H(P)$ is the closure
of the set of all the periodic points whose orbit is homoclinically related to the orbit of $P$.
It also coincides with the closure of all the transverse intersection points of the stable and
the unstable manifolds of the orbit of $P$.
A homoclinic class $H(P)$ is \emph{non-trivial} if it is not
reduced to a periodic orbit.

More generally, if $U$ is an open subset of $M$, we say that two periodic orbits contained in $U$ are
\emph{homoclinically related in $U$} if they are homoclinically related by a pair of heteroclinic points whose orbits stay in $U$.
If $P$ is a periodic point whose orbit is contained in $U$,
we define the \emph{relative homoclinic class $H(P,U)$ of $P$ in $U$} as
the closure of the set of periodic points that are homoclinically related in $U$ to the orbit of $P$.
As above, it coincides with the closure of the set of transverse intersection points between the stable
and the unstable manifolds of $P$ whose orbits stay in $U$. The relative homoclinic class
$H(P,U)$ contains periodic orbits that are related in $U$ to the orbit of $P$ and that are
arbitrarily close to $H(P,U)$ for the Hausdorff distance.\bigskip

We now discuss the link between chain-transitive sets and homoclinic classes
for generic diffeomorphisms of $M$: in this paper, we will consider the set $\diff^1(M)$
of $C^1$-diffeomorphisms of $M$ endowed with the $C^1$-topology; this is a Baire space.
Bonatti and Crovisier established~\cite{BC} a connecting lemma for pseudo-orbits. It implies that
if $f$ belongs to some residual set of $\diff^1(M)$, then,
any chain-recurrence class that contains a periodic orbit is a homoclinic class.
The arguments may be easily localized (see~\cite[theorem 2.2]{crovisier} for
a local version of the connecting lemma for pseudo-orbits).
We thus obtain the following restatement of \cite[remarque 1.10]{BC}.
\begin{theorem}\label{t.rec}
There exists a dense G$_\delta$ subset $\cG_{rec}\subset \diff^1(M)$ such that for any
diffeomorphism $f\in \cG_{rec}$, for any chain-transitive set $K\subset M$ that
contains a periodic point $P$ and for any neighborhood $U$ of $K$, the relative homoclinic
class of $P$ in $U$ contains $K$.
\end{theorem}

Furthermore we proved in~\cite{crovisier} a weak shadowing lemma implying
that the chain-transitive sets that do not contain any periodic point,
(the ``\emph{aperiodic} chain-transitive sets") can also be detected by the periodic orbits:
\begin{theorem}\label{t.shadow}
There exists a dense G$_\delta$ subset $\cG_{shadow}\subset \diff^1(M)$ such that for any
diffeomorphism $f\in \cG_{shadow}$ and any chain-transitive set $K\subset M$,
there exist periodic orbits of $f$ that are arbitrarily close to $K$ for the Hausdorff distance.
\end{theorem}

We end these preliminaries by recalling some standard definitions on the linear subbundles preserved
by a diffeomorphism $f$.

A linear subbundle $E$ of the tangent bundle $T_{K}M$ over an invariant compact set $K$
is \emph{uniformly contracted} by $f$ if there exists an integer $N\geq 1$
such that for any point $x\in K$ and any unit vectors $v\in E_x$ we have
$$\|D_xf^N.v\|\leq \frac{1}{2}.$$
If $E$ is uniformly contracted by $f^{-1}$, we say that $E$ is \emph{uniformly expanded} by $f$.

The tangent bundle on an invariant compact set $K$ has a \emph{dominated splitting}
$T_KM=E\oplus F$ if it decomposes as the sum of two (non-trivial)
linear subbundles $E$ and $F$ and if there exists an integer $N\geq 1$ such that
for any point $x\in K$ and any unit vectors $u\in E_x$ and $v\in F_x$, we have
$$2\|D_xf^N.u\|\leq \|D_xf^N.v\|.$$
An invariant compact set $K$ is \emph{partially hyperbolic} if its tangent bundle
decomposes as the sum of three linear subbundles $E^{ss}\oplus E^c\oplus E^{uu}$
(at most one of the two extremal bundles is non-trivial) so that:
\begin{itemize}
\item If $E^{ss}$ is non-trivial, the splitting $E^{ss}\oplus (E^c\oplus E^{uu})$ is dominated.\\
If $E^{uu}$ is non-trivial, the splitting $(E^{ss}\oplus E^c)\oplus E^{uu}$ is dominated.
\item The bundle $E^{ss}$ is uniformly contracted and the bundle $E^{uu}$ is
uniformly expanded.
\end{itemize}
The bundles $E^{ss}$, $E^{c}$ and $E^{uu}$ are called respectively the \emph{strong stable, central} and
\emph{strong unstable bundles}. 
In the case $E^c$ is trivial, the set $K$ is \emph{hyperbolic}.\\
The partially hyperbolic structure extends to any
invariant compact set contained in a small neighborhood of $K$.
(We refer to~\cite{bonatti-diaz-viana} for a survey on the properties satisfied by partially hyperbolic
sets.)\\
It is known (by M. Brin and Y. Pesin \cite{brin-pesin} or by Hirsch, Pugh and Shub \cite{HPS}),
that at each point $x\in K$ there exists a $C^1$-manifold $W^{ss}(x)$ called the \emph{strong stable manifold} that is tangent to $E^{ss}_{x}$ and characterized by
$$W^{ss}(x)=\left\{y\in M, \; \text{There exists $C>0$ s.t. }\forall n\geq 1, v\in E^{c}_{x},
\;\frac{\dd(f^n(x), f^n(y))}{\|D_xf^n.v\|}\leq C. 2^{-\frac n N}.\right\}.$$
The family $\left(W^{ss}(x)\right)_{x\in K}$ is invariant:
$f(W^{ss}(x))=W^{ss}(f(x))$ for each $x\in K$.
For some small constant $\rho>0$, one can define the \emph{local strong stable manifold} $W^{ss}_{loc}(x)$
at $x\in K$ as the ball centered at $x$ of radius $\rho$ in $W^{ss}(x)$. We get a continuous family
$\left(W^{ss}_{loc}(x)\right)_{x\in K}$ of $C^1$-embedded disks.\bigskip

Let $E$ be an invariant one-dimensional linear subbundle on an invariant compact set $K$.
The diffeomorphism $f$ acts on the bundle of orientations of $E$ (which is also the unit bundle
associated to $E$).  One says that \emph{$f$ preserves an orientation of $E$} if there exists a continuous section
of the bundle of orientations of $E$ over $K$, that is invariant by the action of $f$. In this case, the bundle
$E$ can be identified to the trivial bundle $K\times \RR$. 

The property that $f$ preserves an orientation of the central bundle of a partially hyperbolic set
is open:
\begin{proposition}\label{p.orientation}
Let $f$ be a diffeomorphism and $K$ be a partially hyperbolic set with a one-dimensional
central bundle. If $f$ preserves an orientation on the central bundle, then $f$ preserves
the orientation on the central bundle of any invariant compact set contained in a small
neighborhood of $K$.
\end{proposition}
\begin{proof}
If $f$ preserves an orientation on the central bundle $E^c_K$ over $K$,
then, the bundle of orientations on $E^c_K$ is the union of two disjoint sections
that are preserved by the action of $f$. Let $\La$ be an invariant compact set contained
in a small neighborhood of $K$ (one can assume that $K$ is contained in $\La$).
Then $\La$ is also partially hyperbolic with a one-dimensional central bundle $E^c_\La$.
Moreover, the unit bundles associated to $E^c_\La$ and $E^c_K$ are close (for the Hausdorff
topology). Hence, $f$ also preserves a decomposition of the unit bundle of $E^c_\La$
in two invariant compact sets.
\end{proof}

%%%%%%%%%%%%%%%%%%%%%%%%%%%%%%%%%%%%%%%%%%%%%%%%%%%%%%%%%%%%%%%%%%%%
%%%%%%%%%%%%%%%%%%%%%%%%%%%%%%%%%%%%%%%%%%%%%%%%%%%%%%%%%%%%%%%%%%%%
\section{Abstract central models: definition and properties}\label{s.abstract}
%%%%%%%%%%%%%%%%%%%%%%%%%%%%%%%%%%%%%%%%%%%%%%%%%%%%%%%%%%%%%%%%%%
\subsection{Definition}
We first define in an abstract way a model for the central dynamics.
\begin{definition}
A \emph{central model} is a pair $(\hat K,\hat f)$ where
\begin{itemize}
\item $\hat K$ is a compact metric space (called the \emph{base} of the central model),
\item $\hat f$ is a continuous map from $\hat K\times [0,1]$ into $\hat K\times [0,+\infty)$,
\end{itemize}
and such that
\begin{itemize}
\item $\hat f(\hat K\times \{0\})=\hat K\times \{0\}$,
\item $\hat f$ is a local homeomorphism in a neighborhood of $\hat K\times \{0\}$: there exists a continuous map
$\hat g\colon \hat K\times [0,1]\to \hat K\times [0,+\infty)$ such that
$\hat f\circ \hat g$ and $\hat g\circ \hat f$ are the identity maps
on $\hat g^{-1}(\hat K\times [0,1])$ and $\hat f^{-1}(\hat K\times [0,1])$ respectively,
\item $\hat f$ is a skew product: there exists two maps $\hat f_1\colon \hat K\to \hat K$ and $\hat f_2\colon \hat K\times [0,1]\to [0+\infty)$
respectively such that for any $(x,t)\in \hat K\times [0,1]$, one has
$$ \hat f (x,t)= (\hat f_1(x),\hat f_2(x,t)).$$
\end{itemize}
\end{definition}
Since $\hat f$ preserves the zero section $\hat K\times \{0\}$, it induces a homeomorphism
on the base $\hat K$, also given by the map $\hat f_1$.
The base of the central model is called \emph{minimal} or \emph{chain-transitive} if the dynamics of $\hat f_1$ on $\hat K$ is respectively minimal
or chain-transitive. Since $\hat f$ in general does not preserves $\hat K\times [0,1]$, the dynamics outside
$\hat K\times \{0\}$ is only partially defined. This is however the dynamics we are aimed to analyze here.

%%%%%%%%%%%%%%%%%%%%%%%%%%%%%%%%%%%%%%%%%%%%%%%%%%%%%%%%%%%%%%
\subsection{Existence of trapping strips}
A key assumption for us is the non-existence of chain-recurrent segments in the fibers of the product $\hat K\times [0,1]$.

\begin{definition}
A central model $(\hat K,\hat f)$ has a \emph{chain-recurrent central segment} if it contains
a (non-trivial) segment
%\footnote{One could weaken this definition and also consider a non-trivial segment of the form $\{x\}\times [a,b]$.}
$I=\{x\}\times [0,a]$
contained in a chain-transitive set of $\hat f$.
\end{definition}

\begin{remark}
The fact that the dynamics is only partially defined has no influence on the definition
of a chain-transitive set (see the section~\ref{s.notation}):
in particular, $I$ is contained in the maximal invariant set in $\hat K\times [0,1]$ so that all the iterates
of $I$ are contained in $\hat K\times [0,1]$.
\end{remark}

The subsets of $\hat K\times [0,+\infty)$ we consider will often have the following
geometrical property:

\begin{definition}
A subset $S$ of a product $\hat K\times [0,+\infty)$ is a \emph{strip}
if for any point $x\in \hat K$, the intersection $S\cap \left(\{x\}\times [0,+\infty)\right)$ is an interval.
\end{definition}

We now state the main result of this section:\footnote{From this result, one can consider the dynamics in restriction to the maximal invariant set 
$\Lambda$ in $S$. This set is a strip, but when the base $\hat K$ carries a minimal dynamics,
one can also show that it is bounded by the graph of a continuous map $\varphi\colon \hat K\to [0,1]$, i.e.:
$\La=\{(x,t),\; 0\leq t\leq \varphi(x)\}$.}

\begin{proposition}\label{p.strip}
Let $(\hat K,\hat f)$ be a central model with a chain-transitive base.
Then, the two following properties are equivalent:
\begin{itemize}
\item There is no chain-recurrent central segment.
\item There exist some strips $S$ in $\hat K\times [0,1]$
that are arbitrarily small neighborhoods of $\hat K\times \{0\}$ and
that are trapping regions for $\hat f$ or for  $\hat f^{-1}$:
either $\hat f(\closure(S))\subset \interior(S)$ or $\hat f^{-1}(\closure(S))\subset \interior(S)$.
\end{itemize}
\end{proposition}

This proposition implies that the existence of a chain-recurrent central segment is a local
property on the dynamics in a neighborhood of $K\times \{0\}$.

An open strip $S\subset \hat K\times [0,1)$
satisfying $\hat f(\closure(S))\subset S$ will be called a \emph{trapping strip}.

%%%%%%%%%%%%%%%%%%%%%%%%%%%%%%%%%%%%%%%%%%%%%%%%%%%%%%%%%%%%%%%%%%%%%
\subsection{Proof of proposition~\ref{p.strip}}
Let $\dd_{\hat K}$ be a distance on $\hat K$ and define a distance $\dd$
on $\hat K\times [0,+\infty)$ by
$$\dd((x,s),(y,t))=\max(\dd_{\hat K}(x,y), |t-s|).$$
For $\varepsilon>0$,
one defines the $\varepsilon$-pseudo-orbits of $\hat f$ as the sequences
$(z_0,\dots,z_n)$ in $\hat K\times [0,1]$ such that
$d(f(z_k),z_{k+1})<\varepsilon$ for each $k\in\{0,\dots,n-1\}$.

This allows us to introduce the ``chain-stable" and the ``chain-unstable" sets:
for any invariant subset $\La$ of $\hat K\times [0,1]$ and any $\varepsilon>0$,
we consider:
$$pW^s_\varepsilon(\La)=\{z\in K\times [0,1], \text{ there is a $\varepsilon$-pseudo-orbit }
(z_0,\dots,z_n) \text{ with } z_0=z \text{ and } z_n\in \La\},$$
$$pW^u_\varepsilon(\La)=\{z\in K\times [0,1], \text{ there is a $\varepsilon$-pseudo-orbit }
(z_0,\dots,z_n) \text{ with } z_0\in \La \text{ and } z_n=z\}.$$
One then sets
$$pW^s(\La)=\bigcap_{\varepsilon>0} pW^s_\varepsilon (\La),\quad\text{and} \quad
pW^u(\La)=\bigcap_{\varepsilon>0} pW^u_\varepsilon (\La).$$

One now states three lemmas satisfied by the central models:
\begin{lemma}\label{l.trapping}
For any $\varepsilon>0$, the sets $V^s=pW^s_\varepsilon(\La)$ and $V^u=pW^u_\varepsilon(\La)$
are neighborhoods of $\La$ in $\hat K\times [0,1]$ which satisfy the following properties.\\
$V^s$ is a trapping region for $\hat f^{-1}$: we have
$\hat f^{-1}(\closure(V^s))\cap \left(\hat K\times [0,1]\right) \subset \interior(V^s)$.\\
$V^u$ is a trapping region for $\hat f$: we have
$\hat f(\closure(V^u))\cap \left(\hat K\times [0,1]\right) \subset \interior(V^u)$.
\end{lemma}
\begin{proof}
By definition, the set $V^u=pW^u_\varepsilon(\La)$ contains the $\varepsilon$-neighborhood of $\La$
in $\hat K\times [0,1]$.
Moreover, the $\varepsilon$-neighborhood of $\hat f(V^u)\cap \left(\hat K\times [0,1]\right) $ in $\hat K\times [0,1]$
is contained in $V^u$, showing that
$\hat f(\closure(V^u))\cap \left(\hat K\times [0,1]\right) \subset \interior(V^u)$ for the induced topology on $\hat K\times [0,1]$.

The proof for $V^s=pW^s_\varepsilon(\La)$ is similar.
\end{proof}

\begin{lemma}\label{l.strip}
For any $\varepsilon>0$, the sets $pW^s_\varepsilon(\hat K\times \{0\})$ and
$pW^u_\varepsilon(\hat K\times \{0\})$ are strips of $\hat K\times [0,+\infty)$.

As a consequence, the sets $pW^s(\hat K\times \{0\})$ and
$pW^u(\hat K\times \{0\})$ also are strips of $\hat K\times [0,+\infty)$.
\end{lemma}
\begin{proof}
Let $z=(x,t)$ be a point in $pW^u_\varepsilon(\hat K\times \{0\})$.
There exists a $\varepsilon$-pseudo-orbit $(z_{0},\dots,z_{n})$
such that $z_{0}$ belongs to $\hat K\times \{0\}$ and $z_{n}=z$.

Let us write $z_{k}$ as a pair $(x_{k},t_{k})$ and define
$I_{k}$ as the segment $\{x_{k}\}\times [0,t_{k}]$. Since $z_{k+1}$
is at distance less than $\varepsilon$ from $\hat f(z_{k})$ for each $k$,
one deduces from our choice of the distance $\dd$ that $I_{k+1}$
is contained in the $\varepsilon$-neighborhood of $\hat f(I_{k})$.
The segment $I_{0}=\{z_{0}\}$ is contained in $\hat K\times \{0\}$. Hence, one gets
inductively that $I_{k}$ is included in $pW^u_{\varepsilon}(\hat K\times \{0\})$.

In particular $pW^u_{\varepsilon}(\hat K\times \{0\})$ contains the segment
$\{x\}\times [0,t]=I_n$. This shows that $pW^u_{\varepsilon}(\hat K\times \{0\})$ is a strip.
The proof for $pW^s_{\varepsilon}(\hat K\times \{0\})$ is similar.
\end{proof}

\begin{lemma}\label{l.chain-unstable}
Let us assume that $(\hat K,\hat f)$ has a chain-transitive base and has no chain-recurrent central segment.
Then, the chain-unstable set of $\hat K\times \{0\}$
satisfies one of the following properties:
\begin{itemize}
\item $pW^u(\hat K\times\{0\})=\hat K\times\{0\}$;
\item $\hat K\times [0,\delta]\subset pW^u(\hat K\times\{0\})$ for some $\delta>0$;
\end{itemize}
The same result holds for the chain-stable set of $\hat K\times \{0\}$.
\end{lemma}
\begin{proof}
Let us assume by contradiction that these two properties are not satisfied:
since $pW^u(\hat K\times \{0\})$ is a strip,
one gets (from the first property) a point $x\in \hat K$ and $\eta>0$
such that $\{x\}\times [0,\eta]$ is contained in $pW^u(\hat K\times\{0\})$;
one also gets (from the second one) a point $y\in \hat K$
such that $pW^u(\hat K\times\{0\}) \cap \left(\{y\}\times [0,1]\right)$
is equal to $\{y\}\times \{0\}$.

We fix any constant $\varepsilon_{0}>0$. By our choice of $y$,
for $\varepsilon\in (0,\varepsilon_0)$ small enough,
the strip $pW^u_{\varepsilon}(\hat K\times\{0\})$ does not contain the whole segment
$\{y\}\times [0,\varepsilon_{0}]$.

Using that $\hat K$ is chain-transitive, there exists a $\varepsilon$-pseudo-orbit
$(x_{0},\dots,x_{n})$ in $\hat K$ (for the distance $\dd_{\hat K}$
on $\hat K$) such that $x_{0}=x$ and $x_{n}=y$.
One defines inductively a sequence of segments $I_{k}=\{x_{k}\}\times [0,t_{k}]$
in the following way: $I_{0}$ coincides with $\{x\}\times [0,\eta]$;
if $I_{k}$ has been defined, one denotes by $\{x'_{k}\}\times [0,t'_{k}]$
its image by $\hat f$ and one sets $t_{k+1}=\min(t'_{k},1)$.
In particular, for any $k$, the segment $I_{k+1}$ is both contained in
$\hat K\times [0,1]$ and in the $\varepsilon$-neighborhood of
$\hat f(I_{k})$. Since $I_{0}$ is contained in $pW^u(\hat K\times \{0\})$,
one deduces that $I_{k}$ is contained in $pW^u_{\varepsilon}(\hat K\times \{0\})$
for each $k$. In particular, $I_{n}$ can not contain the segment
$\{y\}\times [0,\varepsilon_{0}]$ so that $t_{n}<\varepsilon_{0}$.

Let $k\in \{0,\dots,n-1\}$ be the smallest integer such that for any $j\in \{k+1,\dots,n\}$,
the segment $I_j$ has a length smaller than $1$ (in other words $t_j<1$).
By this definition, one has $t'_{j-1}=t_j$ for any $j>k$ so that $\hat f(I_{j-1})$ is contained in
the $\varepsilon$-neighborhood of $I_j$. Since $\varepsilon<\varepsilon_0$ and since
$I_n$ is contained in the $\varepsilon_0$-neighborhood of $\hat K\times \{0\}$, one deduces
that $I_k$ and all the $I_j$, for $j>k$ are contained in $pW^s_{\varepsilon_0}(\hat K\times \{0\})$.
From the argument above, $I_k$ is also contained in $pW^u_{\varepsilon_0}(\hat K\times \{0\})$.

By choosing another constant $\varepsilon_0$, one obtains in the same way an interval $I_k=I(\varepsilon_0)$.
Note that by definition, either $k=0$ and $I_k(\varepsilon_0)$ has length $\eta$,
or $k>0$ and $I_k(\varepsilon_0)$ has length $1$.
In any case, when $\varepsilon_0$ goes to $0$, one can extract from the $\left(I_{k}(\varepsilon_0)\right)$
an interval $I=\{z\}\times [0,\eta]$ which is in the chain-recurrence class of $\hat K\times
\{0\}$. In particular, it is a chain-recurrent central segment and
this contradicts our assumption on the central model $(\hat K, \hat f)$.
\end{proof}

We now finish the proof of the proposition \ref{p.strip}.
\begin{proof}[Proof of proposition \ref{p.strip}]
If there exists a chain-recurrent central segment $I$, then, any open strip $S$ which is a trapping region
for $\hat f$ or $\hat f^{-1}$ should contain $I$ and so can not be contained in an arbitrarily small neighborhood of $\hat K\times \{0\}$.

We thus assume conversely that there is no chain-recurrent central segment:
in this case, for any $\delta>0$ the strip $\hat K\times [0,\delta]$ can not be both contained in
$pW^s(\hat K\times\{0\})$ and in $pW^u(\hat K\times\{0\})$. By lemma \ref{l.chain-unstable}, one deduces that
$pW^s(\hat K\times\{0\})$ or $pW^u(\hat K\times\{0\})$ is equal to $\hat K\times\{0\}$.
We will suppose that we are in the second case (in the first case, the proof is similar).

Since $pW^u_{\varepsilon}(\hat K\times\{0\})$ is decreasing
towards $pW^s(\hat K\times\{0\})=\hat K\times\{0\}$
when $\varepsilon$ goes to $0$, one deduces that for $\varepsilon>0$ small enough, the set
$S= pW^u_{\varepsilon}(\hat K\times\{0\})$ is an arbitrarily small neighborhood of $\hat K\times\{0\}$.
By lemma \ref{l.strip}, $S$ is a strip and by lemma \ref{l.trapping}, we have
$$\hat f(\closure(S))\cap \left(\hat K\times [0,1]\right) \subset \interior(S) \text{ in } \hat K\times [0,1].$$
Since $S$ and $\hat f(S)$ are small neighborhood of $\hat K\times\{0\}$, this also gives
the required property: $\hat f(\closure(S))\subset \interior(S)$ for the topology of $\hat K\times [0,+\infty)$.
\end{proof}

%%%%%%%%%%%%%%%%%%%%%%%%%%%%%%%%%%%%%%%%%%%%%%%%%%%%%%%%%%%%%%%%%
%%%%%%%%%%%%%%%%%%%%%%%%%%%%%%%%%%%%%%%%%%%%%%%%%%%%%%%%%%%%%%%%%
\section{Central models for partially hyperbolic dynamics}\label{s.application}
%%%%%%%%%%%%%%%%%%%%%%%%%%%%%%%%%%%%%%%%%%%%%%%%%%%%%%%%%%%%%%%%%
\subsection{Definition of central models for partially hyperbolic dynamics}
We here explain how central models are related to partially hyperbolic dynamics.
\begin{definition}
Let $f$ be a diffeomorphism of a manifold $M$ and $K$ be a partially hyperbolic invariant compact set for $f$
having a one-dimensional central bundle. A central model $(\hat K,\hat f)$
is a \emph{central model for $(K,f)$} if there exists a continuous map
$\pi\colon \hat K\times [0,+\infty)\to M$ such that:
\begin{itemize}
\item $\pi$ semi-conjugates $\hat f$ and $f$: we have $f\circ \pi=\pi\circ \hat f$ on
$\hat K\times [0,1]$,
\item $\pi(\hat K\times \{0\})=K$,
\item The collection of maps $t\mapsto \pi(\hat x,t)$ is a continuous family of $C^1$-embeddings
of $[0,+\infty)$ into $M$, parameterized by $\hat x\in \hat K$.
\item For any $\hat x\in \hat K$, the curve $\pi(\hat x,[0,+\infty))$ is tangent at the point $x=\pi(\hat x,0)$
of $K$ to the central bundle.
\end{itemize}
\end{definition}

One can associate a central model to any partially hyperbolic set.

\begin{proposition}\label{p.existence}
Let $f$ be a diffeomorphism and $K$ be a partially hyperbolic invariant compact set having a one-dimensional central bundle.
Then, there exists a central model for $(K,f)$.

Moreover, one can choose the central model $(\hat K,\hat f)$ such that if $K$ is chain-transitive
(resp. minimal), then the base of the central model is also chain-transitive (resp. minimal).
\end{proposition}
This will be proven in the next section.
Following the dichotomy of proposition~\ref{p.strip},
we then discuss in sections~\ref{ss.segment} and~\ref{ss.strip}
the cases the central model has a chain-recurrent central segment or attracting strips.

%%%%%%%%%%%%%%%%%%%%%%%%%%%%%%%%%%%%%%%%%%%%%%%%%%%%%%%%%%%%%%%%%%%%%%%%%%%%%%%%%%%%%%%%%%%%%%%%%%%%%%%%%%
\subsection{Construction of central models}
We here prove proposition~\ref{p.existence} giving the existence of central models for partially hyperbolic
sets. This is a consequence of the plaque family theorem of Hirsch, Pugh and Shub.

Let us consider a compact subset $K$ of a manifold $M$ and a linear subbundle $E_K$ of
the tangent bundle $T_KM$ above $K$. At each point $x\in K$ and for $r>0$, we denote by $E_x(r)$ the open ball
of radius $r$ of the vector space $E_x$ and by $E_K(r)$ the collection of all the balls $E_x(r)$ for $x\in K$.

A \emph{plaque family tangent to $E$} is a continuous map $\cD$ from $E_K(1)$ into $M$ such that:
\begin{itemize}
\item For each point $x\in K$, the induced map $\cD_x\colon E_x(1)\to M$ is a $C^1$-embedding.
\item $\cD_x(0)=x$ and $E_x$ is the tangent space at $x$ to the image of $\cD_x$.
\item $(\cD_x)_{x\in K}$ is a continuous family of $C^1$-embeddings of the disks $(E_x(1))_{x\in K}$.
\end{itemize}

When $K$ and $E_K$ are invariant by a diffeomorphism $f$ and its tangent map $D_Kf$, one says that the plaque family is \emph{locally invariant}
if there exists $r>0$ satisfying the following property:
for each $x\in K$, the image by $f$ (resp. $f^{-1}$)
of the embedded disk $\cD_x(E_x(r))$ of radius $r$ at $x$
is contained in the embedded unit disk $\cD_{f(x)}(E_{f(x)}(1))$ at $f(x)$
(resp. $\cD_{f^{-1}(x)}(E_{f^{-1}(x)}(1))$ at $f^{-1}(x)$).

We now state the plaque family theorem by Hirsch, Pugh and Shub (see~\cite{HPS} theorem 5.5).
\begin{plaque}
Let $K$ be an invariant compact set for a diffeomorphism whose tangent space decomposes
into a dominated splitting $T_KM=E\oplus F$. Then, there exists a locally invariant plaque family that
is tangent to $E$.
\end{plaque}

\begin{remark}\label{r.plaquecentral}
If the tangent space of $K$ has a dominated splitting into three bundles $T_KM=E\oplus F\oplus G$,
then, the theorem gives the existence of a plaque family tangent to $E\oplus F$ and another one
tangent to $F\oplus G$. At each point of $K$, one can take the intersections between the plaques
of these two families and obtain a locally invariant plaque family that is tangent to $F$.
\end{remark}

\begin{proof}[End of the proof of proposition~\ref{p.existence}]
Let us consider an invariant compact set $K$ which is partially hyperbolic (we have a splitting
$T_KM=E^{ss}\oplus E^c \oplus E^{uu}$) and whose central bundle $E^c$ is one dimensional.
By using the plaque family theorem (and the remark~\ref{r.plaquecentral}), we get a locally invariant
plaque family $\cD$ tangent to $E^c$.

One will first assume that the bundle $E^c$ is orientable and that $f$ preserves an orientation
of $E^c$. Since the bundle $E^c$ is one-dimensional, it is a trivial bundle
and $E^c(1)$ can be identified with the product $K\times \RR$ so that the map $\cD$ sends $K\times \{0\}$
onto $K$. Since the plaque family is locally invariant, the dynamics of $f$ in $M$
lifts as a continuous dynamical system $\hat f\colon U\to K\times \RR$ defined on a neighborhood
$U$ of $K\times \{0\}$ in $K\times \RR$:
we have $\cD\circ \hat f=f \circ \cD$ and $\hat f$ is a skew product on $U$.
Up to changing the trivialization of $E^c$,
one may assume that the canonical orientation on $K\times \RR$ is preserved by $\hat f$.
Multiplying by a constant along the second coordinate of $K\times \RR$, one
can assume that $\hat f$ is defined on $K\times [-1,1]$. Since $\hat f$ preserves the orientation,
it sends $K\times [0,1]$ into $K\times [0,+\infty)$.
Since the same arguments apply to $f^{-1}$, one gets that $\hat f$ is a local homeomorphism.
Hence the map $\cD\colon K\times [0,+\infty)\to M$ is the projection $\pi$ of the central model
$(K,\hat f)$ for $(K,f)$,
finishing the proof of the proposition in the orientable case.
(The minimality and the chain-transitivity are preserved by the construction.)

When $E^c$ has no orientation preserved by $f$,
one can consider the two-fold orientation covering $\tilde E^c$ of $E^c$
and the lifted map $\tilde \cD$ from $\tilde E^c(1)$ to $M$.
The dynamics on $M$ again lifts as a continuous dynamical system $\hat f$ defined in a neighborhood $U$
of the zero section of $\tilde E^c$. Moreover, considering any
orientation of the bundle $\tilde E^c$, one can require that the map $\hat f$ preserves this
orientation (as a consequence $\hat f$ is uniquely defined in $U$ by the choice of an orientation).
Denoting by $\tilde K$ the base (i.e. the zero section) of $\tilde E^c$,
one identifies $\tilde E^c(1)$ with the product $\tilde K \times \RR$ and as above,
the map $\tilde \cD\colon \tilde K\times [0,+\infty)\to M$ shows that $(\tilde K,\hat f)$
is a central model for $(K,f)$.
The map $\tilde \cD\colon \tilde K\times \{0\}\to K$ can be viewed as the bundle of orientations
on $E^c_{K}$. Hence, the last assertion of the proposition is a consequence of the next lemma.
\end {proof}

\begin{lemma}\label{l.orientation} Let $K$ be an invariant compact set and
$E$ be an invariant one-dimensional linear subbundle of the tangent bundle on $K$.
Let $\hat f$ be the lift of $f$ on $\pi\colon \hat K\to K$, the bundle of orientations of $E$ over $K$.

If $K$ is chain-transitive (resp. minimal) and if $f$ does not preserve any orientation on $E$, then,
$(\hat K,\hat f)$ is also chain-transitive (resp. minimal).
\end{lemma}
\begin{proof}
Any point $x\in K$ has two preimages $\hat x^-$ and $\hat x^+$ by $\pi$ and there exists
a homeomorphism $\sigma\colon \hat K\to \hat K$ which exchanges $\hat x^-$ and $\hat x^+$ for each
$x\in K$. The assumption that $E$ has no orientation preserved by $f$ can be restated as:
$\hat K$ is not the disjoint union of two invariant compact sets $A$ and $A'=\sigma(A)$.

Let us now assume that $(K,f)$ is chain-transitive: for any two points $x,y\in K$, we have
$x\dashv y$ (see section~\ref{s.notation}). By lifting the pseudo-orbits to $\hat K$, we get that one of the two properties is satisfied:
\begin{itemize}
\item $\hat x^-\dashv \hat y^-$ and $\hat x^+\dashv \hat y^+$ ;
\item $\hat x^-\dashv \hat y^+$ and $\hat x^+\dashv \hat y^-$ .
\end{itemize}
One deduces that any chain recurrence class of $\hat K$ projects down by $\pi$ onto the whole $K$.
Hence, if $\hat K$ is not chain-transitive, it is the disjoint union of two chain transitive classes
$A$ and $A'=\sigma(A)$ giving a contradiction.

We end the proof by assuming that $(K,f)$ is minimal. Any minimal set $A\subset \hat K$ projects down
by $\pi$ onto the minimal set $K$. Moreover $A'=\sigma(A)$ is also a minimal set. Hence, either
$A=A'=\hat K$ showing the minimality of $\hat K$ or $A\cap A'=\emptyset$ giving a contradiction.
\end{proof}

Let us remark that if $x$ is a point in $K$, then, any point $\hat x\in \pi^{-1}(x)$ in
the central model is associated to a \emph{half central curve}: the parameterized $C^1$-curve $\pi(\hat x, [0,+\infty))$.
It defines an orientation of $E^c_x$.

The construction of $(\hat K,\hat f)$ in the proof of proposition \ref{p.existence} gives more information:

\begin{addendum}\label{a.model}
One can build the central model $(\hat K,\hat f)$ of proposition~\ref{p.existence}
with the following additional properties:
\begin{enumerate}
\item\label{i.model1}The curves $\pi(\{\hat x\}\times [0,+\infty))$, $x\in \hat K$,
are contained in the plaques of a prescribed plaque family tangent to $E^c$.
In particular, by remark~\ref{r.plaquecentral}
one may assume that each of these curves is the intersection of
a plaque tangent to $E^{ss}\oplus E^c$ with a plaque tangent to $E^c\oplus E^{uu}$.
\item\label{i.model2}  The projection $\pi(\hat K\times [0,+\infty))$ of the central model can be chosen
in an arbitrarily small neighborhood of $K$. As a consequence, one may assume that the partially
hyperbolic structure of $K$ extends to the projection
$\Lambda=\pi(\hat \Lambda)$ of the maximal invariant set $\hat\Lambda=\cap_{n\in \ZZ} \hat f^n(\hat K\times [0,1])$;
in the case $f$ preserves an orientation of $E^c_K$, it can also preserve an orientation of $E^c$
on $\Lambda$.
\item\label{i.model3} If $f$ preserves an orientation of $E^c_K$, then $\pi$ is a bijection between
$\hat K$ and $K$. By considering the opposite orientation on $E^c_K$, one can build
another central model $(\hat K',\hat f')$ for $(K,f)$ with a projection $\pi'$. Hence, any point $x\in K$ is contained
in a central curve which is the union of the two half central curves $\pi(\{\hat x\}\times [0,+\infty))$ and $\pi'(\{\hat x\}\times [0,+\infty))$
at $x$ given by the two models for $(K,f)$.
\item\label{i.model4} If $f$ does not preserve any orientation of $E^c_K$, then $\pi\colon \hat K\to K$ is two-to-one:
any point $x\in K$ has two preimages $\hat x^-$ and $\hat x^+$ in $\hat K$.
The homeomorphism $\sigma$ of $\hat K$ which exchanges the preimages
$\hat x^-$ and $\hat x^+$ of any point $x\in K$ commutes with $\hat f$.\\
Any point $x\in K$ is contained
in a central curve which is the union of the two half central curves $\pi(\{\hat x^-\}\times [0,+\infty))$ and $\pi(\{\hat x^+\}\times [0,+\infty))$ at $x$
associated to $\hat x^-$ and to $\hat x^+$ (they have different orientation).
\end{enumerate}
\end{addendum}

%%%%%%%%%%%%%%%%%%%%%%%%%%%%%%%%%%%%%%%%%%%%%%%%%%%%%%%%%%%%%%%%%%%%%%%%%%%%%%%%%%%%%%%%%%%%%%%%%%%%%
\subsection{Chain-recurrent central segments}\label{ss.segment}
In this section we consider the existence of a chain-recurrent central segment.
\begin{proposition}\label{p.segment}
For any diffeomorphism $f$ in the dense G$_\delta$ subset $\cG_{rec}\cap\cG_{shadow}$
of $\diff^1(M)$, any chain-transitive compact set $K$
which is partially hyperbolic with a one-dimensional central bundle
satisfies the following property:

If a central model of $(K,f)$ has a chain-recurrent central segment,
then $K$ is contained in (non-trivial) relative homoclinic classes that are included
in arbitrarily small neighborhoods of $K$.
\end{proposition}

One gets this proposition from the following result, which is a restatement of theorem 12 of~\cite{bonatti-gan-wen} (we state it for any chain-transitive set).

\begin{proposition}\label{p.central}
Let $f$ be a diffeomorphism in the dense G$_\delta$ subset $\cG_{rec}\cap\cG_{shadow}$
and $\La$ be a chain-transitive and partially hyperbolic set with a one-dimensional central bundle.

If $\La$ contains a $C^1$ curve $\gamma$ that is not tangent to the strong stable
nor to the strong unstable bundle, then, $\La$ is contained in a (non-trivial) homoclinic class
$H(P)$. More precisely: in any neighborhood of $\La$ there exists a relative homoclinic class
containing $\La$.
\end{proposition}

\begin{figure}[ht]
\begin{center}
\input{central.pstex_t}
\end{center}
\caption{Strong invariant manifolds in the neighborhood of a chain-recurrent central segment.\label{f.central}}
\end{figure}
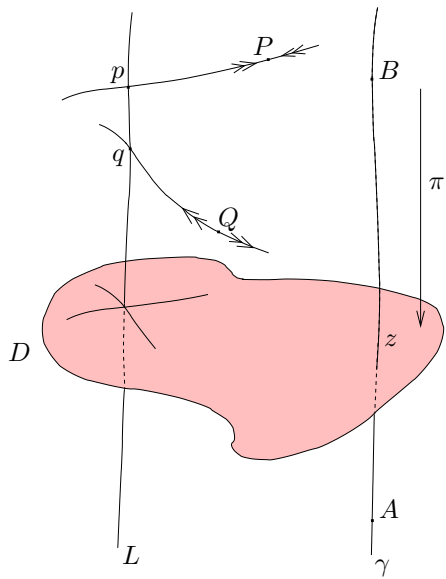

Proposition~\ref{p.central} is a consequence of the results in~\cite{BC} and~\cite{crovisier}.
For completeness we explain the main ideas of the proof, as given by Bonatti, Gan and Wen:

Let us first discuss the case both bundles $E^{ss}$ and $E^{uu}$ are non trivial (see figure~\ref{f.central}).
Let us consider a small segment $[A,B]$ in the curve $\gamma$,
a point $z$ in the interior of this segment,
and a small differentiable disk $D$ centered at $z$ and transverse to $\gamma$. Any
one-dimensional foliation, that contains the curve $\gamma$ as a leaf, defines (by holonomy along the foliation) a submersion $\Pi$ from a
neighborhood of $\gamma$ onto $D$. This foliation may be endowed with an orientation on its leaves
such that $A<B$ on the segment $\gamma$.\\
Let $P$ and $Q$ be two points close to $\gamma$ and whose orbits stay in a small neighborhood of $\Lambda$.
They have local strong stable and strong unstable manifolds. The projections of $W^{ss}_{loc}(P)$ and $W^{uu}_{loc}(Q)$ by $\Pi$
intersect in a unique point of $D$. This shows that there exists a unique leaf $L$ of the foliation which intersects
both $W^{ss}_{loc}(P)$ and $W^{uu}_{loc}(Q)$. Each intersection is reduced to a unique point $p$ and $q\in L$ respectively.
Using the order in the leaf $L$, one denotes $W^{ss}_{loc}(P)<W^{uu}_{loc}(Q)$ if $p<q$ and $W^{ss}_{loc}(P)>W^{uu}_{loc}(Q)$ if $q>p$.
Note that the pair $(p,q)$ varies continuously with $(P,Q)$. If one fixes the point $Q$ close to $z$ and if
$P$ varies inside the segment $[A,B]$, one gets $W^{ss}_{loc}(A)<W^{uu}_{loc}(Q)$ and $W^{ss}_{loc}(B)>W^{uu}_{loc}(Q)$.
Hence, by continuity, there exists a point $P$ in $\gamma$ such that $W^{ss}_{loc}(P)$ and $W^{uu}_{loc}(Q)$ intersect.

The case $E^{uu}$ (or $E^{ss}$) is trivial is simpler: for any point $Q$ close to $z$ and whose orbit stays in a small neighborhood
of $\Lambda$, the (one-codimensional) strong stable manifold of $Q$ intersects $\gamma$. One deduces that the forward orbit
of $Q$ accumulates on a subset of $\Lambda$ in the future.

We thus have proven:
\begin{lemma}\label{l.intermediate}
Let $U$ be a neighborhood of $\Lambda$ and $z$ a point in $\gamma$.
Then, there exists a neighborhood $V$ of $\Lambda$ and $W$ of $z$ such that
for any point $Q\in W$ whose orbit stays in $V$, either the orbit of $Q$ accumulates
(in the past or in the future) on a subset of $\Lambda$, or both bundles $E^{ss}$ and $E^{uu}$
are non-trivial and there exist
\begin{itemize}
\item a point $p^s\in W^{ss}(Q)$ whose orbit stays in $U$ and accumulates on a subset of $\Lambda$ in the past,
\item a point $p^u\in W^{uu}(Q)$ whose orbit stays in $U$ and accumulates on a subset of $\Lambda$ in the future.
\end{itemize}

In particular, any chain-transitive set contained in $V$ and intersecting $W$ is contained in
a larger chain-transitive set $\Lambda'\subset \closure(U)$ that intersects $\Lambda$.
\end{lemma}

Since $f$ belongs to $\cG_{shadow}$, by theorem~\ref{t.shadow} the set $\Lambda$
is approximated for the Hausdorff distance by periodic orbits having
a point close to some point $z\in \gamma$. One deduces that there exist
chain-transitive sets $\Lambda'$ that contain $\Lambda$, that contain some periodic orbits,
and that are contained in an arbitrarily small neighborhood of $\Lambda$.
Since $f$ belongs to $\cG_{rec}$, by theorem~\ref{t.rec}
the set 
$\Lambda'$ can be assumed to be a relative homoclinic class, concluding
the proof of the proposition \ref{p.central}.

Note that the proof gives some interesting additional properties (that we will not use here):
\begin{addendum}\label{a.cycle}
Under the assumptions of proposition~\ref{p.central}, we have moreover:
\begin{itemize}
\item Any periodic orbit which is close to $\Lambda$ for the Hausdorff distance belongs
to the homoclinic class $H(P)$.
\item If both extremal bundles $E^{ss}$ and $E^{uu}$ on $\Lambda$ are non-trivial, then, by arbitrarily small $C^1$-perturbations, one can create strong connections
for the periodic orbits close to $\La$ for the Hausdorff distance:
that is some intersections between the strong stable and
the strong unstable manifolds of the periodic orbit.
\end{itemize}
If the central Lyapunov exponent of any of these periodic orbit is weak, one can create
a heterodimensional cycle.
\end{addendum}

We now prove proposition~\ref{p.segment}.
\begin{proof}[Proof of proposition~\ref{p.segment}]
Let $(\hat K,\hat f)$ be a central model for $(K,f)$
with a chain-recurrent central segment.
By proposition~\ref{p.strip}, the existence of
chain-transitive central segments is a local property of $\hat K$: if $(\hat K,\hat f)$ has a chain-transitive
central segment, then, in any neighborhood of $\hat K\times \{0\}$ in
$\hat K\times [0,1)$, there exists an invariant compact
chain-transitive set $\hat \La$ which contains a chain-transitive central segment $I$.

If one projects $\hat \La$ in the manifold, one obtains a chain-transitive set $\La$ contained
in an arbitrarily small neighborhood of $K$. Hence $\La$ is partially hyperbolic.

The projection of $I$ in the manifold is a $C^1$-curve $\gamma$ which is tangent to the central bundle
at some point of $K\cap \gamma$. Hence, the curve $\gamma$ is not tangent to the strong stable nor to the strong unstable bundles and proposition~\ref{p.central} applies to $\La$:
$\Lambda$ is contained in a relative homoclinic class $H(P, U)$.
Since $\Lambda$ can be chosen arbitrarily close
to $K$, the set $H(P,U)$ is contained in an arbitrarily small neighborhood $U$ of $K$.
\end{proof}

%%%%%%%%%%%%%%%%%%%%%%%%%%%%%%%%%%%%%%%%%%%%%%%%%%%%%%%%%%%%%%%%%%%%%%%%%%%%%%%%%%%%%%%%%%%%%%%%%%%%%%%%
\subsection{Trapping strips}\label{ss.strip}
In this section we discuss the dynamics of central models having a trapping strip.
\begin{proposition}\label{p.diff-strip}
Let $f$ be a Kupka-Smale diffeomorphism and $K$ be a compact set
which is partially hyperbolic with a one-dimensional central bundle and
which contains infinitely many periodic points.
Let $(\hat K,\hat f)$ be a central model for $(K,f)$ that satisfies the properties stated in addendum~\ref{a.model}
and let $\pi\colon \hat K\times [0,+\infty)\to M$ be
the associated projection.

If the central model $(\hat K,\hat f)$ has a trapping strip $S$, then, for any
neighborhoods $\cU$ of $f$ in $\diff^1(M)$ and $U$ of $\pi(S)$, there exist
\begin{itemize}
\item a periodic orbit $\cO$ of $f$ that is the projection by $\pi$ of a periodic orbit
$\hat \cO\subset S$ of $\hat f$,
\item a diffeomorphism $g\in \cU$ that coincides with $f$ outside $U$ and on a neighborhood of $\cO$
\end{itemize}
such that the relative homoclinic class of $\cO$ in $U$ is non-trivial.
\end{proposition}

In the proof of the proposition, the bundles $E^{ss}$ and $E^{uu}$ will
play a different role, due to our assumption on
the existence of a trapping strip $S$ for $f$ (not for $f^{-1}$).
Before, we discuss the local dynamics in the neighborhoods of the periodic points of $K$
(see lemma~\ref{l.periodic} below).

For any point $\hat x\in \hat K$, we denote by
$\sigma_{\hat x}$ the curve in $M$ which is the projection by $\pi$
of the segment $S\cap \left(\{\hat x\}\times [0,1)\right)$.
Since $S$ is a trapping region, we have
the inclusion
$f(\closure(\sigma_{\hat x}))\subset \sigma_{\hat f(\hat x)}$.

One can also consider the maximal invariant set
$\Lambda=\cap_{n\geq 0}\hat f^n(S)$ in $S$. This is a compact
strip and for any point $\hat x\in \hat K$, we denote by
$\gamma_{\hat x}\subset \sigma_{\hat x}$
the curve in $M$ which is the projection by $\pi$
of the segment $\Lambda\cap \left(\{\hat x\}\times [0,1)\right)$.
We have the equality
$f(\gamma_{\hat x})= \gamma_{\hat f(\hat x)}$.

Since $\Lambda$ is compact and contained in the open set $\hat f(S)$,
the family $(\gamma_{\hat x})_{\hat x\in \hat K}$ satisfies a semi-continuity property:
for any $\delta>0$ and for any points $\hat x, \hat x'\in \hat K$ that are close enough,
the curve $\gamma_{\hat x'}$ is contained in the
$\delta$-neighborhood of $\sigma_{\hat x}$.\bigskip

We state some consequences of the properties of addendum \ref{a.model}:
\begin{itemize}
\item There exists a plaque family $(\cD^{cs}_x)_{x\in K}$ tangent to $E^{ss}\oplus E^c$
such that any curve $\gamma_{\hat x}$ is contained in the plaque $\cD^{cs}_{x}$
centered at the point $x=\pi(\hat x)$.
\item $\Lambda$ is partially hyperbolic with a one-dimensional central bundle.
\item Any preimage $\hat x$ by $\pi$ of a periodic point $x\in K$ is periodic.
In particular, the curve $\gamma_{\hat x}$ is periodic by $f$.
One endpoint is $x=\pi(\hat x)$, the other one is a periodic point denoted by $P_{\hat x}$.\\
Since $f$ is Kupka-Smale, $\gamma_{\hat x}$ is the union of finitely many periodic
points and of segments contained in the stable sets of these periodic points.
\end{itemize}
If $z$ is a periodic point contained in a periodic segment $\gamma_{\hat x}$,
it belongs to the partially hyperbolic set $\Lambda$. Moreover,
its local strong stable manifold $W^{ss}_{loc}(z)$ is contained in the plaque $\cD^{cs}_{x}$.
(The local strong stable manifold can be obtained as the unique fixed point
of a graph transform.
Since the plaque family $\cD^{cs}$ is locally invariant, the iterates of a graph
in a plaque of $\cD^{cs}$ under the graph transform converge towards
the local strong stable manifold in the plaques of $\cD^{cs}$).
Since $W^{ss}_{loc}(z)$ is one-codimensional in $\cD^{cs}_{x}$, it separates the plaque
in two components. In particular, one can define
the half plaque $\widehat \cD^{cs}_{\hat x}$ bounded by $W^{ss}_{loc}(x)$ and containing the curve
$\gamma_{\hat x}$.

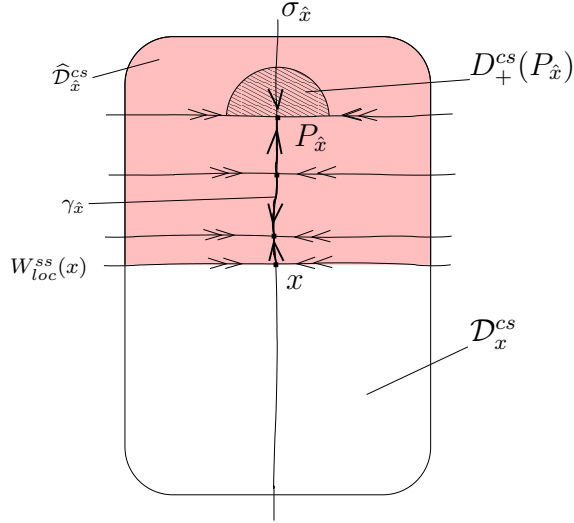
\begin{figure}[ht]
\begin{center}
\input{periodic.pstex_t}
\end{center}
\caption{The center-stable plaque of a periodic point $x\in K$.\label{f.periodic}}
\end{figure}
We describe the dynamics inside the half periodic plaques (see figure~\ref{f.periodic}):
\begin{lemma}\label{l.periodic}
There exists $\eta>0$ such that for any point $\hat x\in \hat K$ which lifts a periodic point
$x\in K$, one has:
\begin{enumerate}
\item The union of the local stable manifolds of the periodic points that belong to
$\gamma_{\hat x}$ contains the $\eta$-neighborhood of $\sigma_{\hat x}$ in the half plaque
$\widehat \cD^{cs}_{\hat x}$.
\item In $\cD^{cs}_{x}$, the local stable manifold of $P_{\hat x}$ contains the half open ball
$D^{cs}_{+}(P_{\hat x})$
centered at $P_{\hat x}$ and of radius $\eta$ which is bounded by the local strong manifold
$W^{ss}_{loc}(P_{\hat x})$ and which is both contained in $\hat \cD_{\hat x}^{cs}$ and disjoint from $\gamma_{\hat x}$.
\end{enumerate}
\end{lemma}
\begin{proof}
The proof is standard:
one first considers a cone field $\cC^{ss}$ around the bundle
$E^{ss}$. This cone field is invariant and its elements are expanded by some iterate $f^{-N}$. Furthermore, if one defines
$\bar \sigma_{x}=f^{-N}(\sigma_{f^{N}(x)})$ we have 
$\closure(\sigma_{\hat x})\subset \bar \sigma_{\hat x}$.

One fixes a small constant $\delta>0$ and
for any periodic point $\hat x$,
one considers the $\delta$-neighborhood $U(\hat x)$
of $\sigma_{\hat x}$ in the half plaque $\widehat \cD^{cs}_{\hat x}$:
any point $z\in U(\hat x)$
can be joined to a neighborhood of $\sigma_{\hat x}$
in $\bar \sigma_{\hat x}$ by a small curve $c$
tangent to $\cC^{ss}$ and contained in $\widehat \cD^{cs}_{\hat x}$.

Let us consider an orbit $(z_0,\dots,z_n)$
contained in the neighborhoods $U(\hat x),\dots,U(\hat f^n(\hat x))$
and such that $z_0$ belongs to the $\eta$-neighborhood of $\sigma_{\hat x}$,
for some small $\eta>0$.
The backward iterate $f^{-n}(c_n)$ of a small curve $c_n$
tangent to $\cC^{ss}$ that joins $z_n$ to
a point $z'_n$ of $\bar \sigma_{\hat f(\hat x)}$
is tangent to $\cC^{ss}$, contained in $U(\hat x)$ and has a bounded length.
This implies (by contraction of the vectors in the cone field $\cC^{ss}$) that the length of $c_n$ is exponentially small when $n$ is large.
Moreover, the point $f^{-n}(z'_n)$ belongs to $\bar \sigma_{\hat x}$ so that
$z'_n$ is arbitrarily close to $\gamma_{\hat x}$ when $n$ is large.
We have thus showed that the forward orbit of any point $z_0$ in the
$\eta$-neighborhood of $\sigma_{\hat x}$ in the half plaque
$\widehat \cD^{cs}_{\hat x}$ accumulates on a subset
$\gamma_{\hat x}\cup \dots f^{\tau-1}(\gamma_{\hat x})$,
where $\tau$ is the period of $\hat x$ by $\hat f$.
One concludes the proof using that the dynamics on $\gamma_{\hat x}$ is Morse-Smale.
\end{proof}

\begin{corollary}\label{c.strong-unstable}
The bundle $E^{uu}$ is non-trivial.
\end{corollary}
\begin{proof}
The proof is done by contradiction: if $E^{uu}$ is trivial, lemma~\ref{l.periodic}
shows that for any periodic point $\hat x\in \hat K$, the associated periodic point $P_{\hat x}$ is
a sink whose basin (the points $y\in M$ such that $d(f^n(P_{\hat x}),f^n(y))\rightarrow 0$ when
$n\rightarrow +\infty$) has a volume uniformly bounded from below.
The basin of distinct iterates of $P_{\hat x}$ are pairwise disjoint and since $M$ has finite
volume, the period of $P_{\hat x}$ is uniformly bounded. Since $f$ is Kupka-Smale,
and $K$ contains infinitely many periodic orbits, the periodic points in $\hat K$ have
arbitrarily large period. The same property holds for the points $P_{\hat x}$ giving a contradiction.
\end{proof}

We now come to the proof of proposition~\ref{p.diff-strip}:
two cases should be considered.

\subsubsection{The non-orientable case}
By addendum~\ref{a.model}.(\ref{i.model3}), any point $x\in K$ has two preimages $\hat x^-$ and $\hat x^+$ by $\pi$.
Moreover, the set $\Gamma_{x}=\gamma_{\hat x^-}\cup \gamma_{\hat x^+}$ is a $C^1$-curve
contained in the plaque $\cD^{cs}_{x}$. By corollary~\ref{c.strong-unstable},
one can consider the local strong unstable manifolds of points in $\Gamma_{x}$:

\begin{lemma}\label{l.return}
Let $x,x'$ be two periodic points in $K$ that are close to each other.
Then, the local strong unstable manifold of any point $y\in \Gamma_{x'}$
intersects transversally the local stable manifold of a periodic point of $\Gamma_{x}$.
\end{lemma}
\begin{proof}
Since, $x$ and $x'$ are close, the local strong unstable manifold
$W^{uu}_{y}$ intersects the plaque $\cD^{cs}_{x}$.
By semi-continuity of the family $(\gamma_{\hat z})_{\hat z\in \hat K}$,
the intersection point belongs to the $\eta$-neighborhood of $\sigma_{\hat x^-}\cup \sigma_{\hat x^+}$
in $\cD_{x}^{cs}$. Hence, by lemma~\ref{l.periodic}, the strong unstable manifold of
$z$ intersects transversally the stable manifold of a periodic point
of $\Gamma_x$.
\end{proof}

The diffeomorphism $f$ is a Kupka-Smale diffeomorphism; since $K$ contains infinitely many periodic points,
it contains periodic points with arbitrarily large period.
In particular, there exists a periodic point $x\in K$ having a return $x'=f^n(x)\neq x$
arbitrarily close to $x$ and lemma~\ref{l.return} may be applied to the pair $(x,f^n(x))$.
One will work with the periodic points $p$ in $\Gamma_x$: each of them may be lifted
by $\pi^{-1}$ in $\{\hat x^-\}\times [0,+\infty)$
or in $\{\hat x^+\}\times [0,+\infty)$ as a periodic point for $\hat f$.
One will show that one of them has a non-trivial relative homoclinic class in a small neighborhood
$U$ of $\pi(S)$.

Let us consider Smale's partial order on the periodic points:
for two periodic points $p$ and $q$, we write $p\geq q$ if the unstable
manifold of $p$ intersects transversally the stable manifold of $q$.
The $\lambda$-lemma implies that this is an order.
One chooses among the (finitely many) periodic point contained in $\Gamma_x$
a point $p$ which is minimal for Smale's order. By lemma~\ref{l.return},
there exists a periodic point $q$ in $\Gamma_{x'}$ such that $p\geq q$.
Similarly, there exists a point $p'$ in $\Gamma_x$ such that $q\geq p'$.
Since $p$ is minimal for Smale's ordering among the
periodic points of $\Gamma_x$, one also has $p'\geq p$; by transitivity
of $\geq$ one thus gets $q\geq p$.
We thus have $p\geq q\geq p$, showing that the two distinct periodic points $p,q$ are homoclinically
related. This implies that their homoclinic classes coincide and are non-trivial.

Note that if one chooses the points $x,x'\in\pi(S)$ close enough, the points $p$ and $q$ are related
by a pair of heteroclinic points whose orbits stay in an arbitrarily small open neighborhood $U$ of $\pi(S)$: 
one defines $p\geq q$ for periodic points whose orbits are contained in $U$
if there exists a transverse intersection point between the unstable manifold of $p$
and the stable manifold of $q$ whose orbit is contained in $U$. Since $U$ is open, by the $\lambda$-lemma,
this defines again an order on the periodic points whose orbits are contained in $U$. 
This proves proposition~\ref{p.diff-strip} in the non-orientable case (we did
not need to perturb $f$).

\subsubsection{The orientable case}\label{sss.orientable}
In this case, one has to consider points in $K$ that are in a twisted position.
This geometry already appeared in~\cite{bonatti-gan-wen} and we describe here
a perturbation result Bonatti, Gan and Wen proved there for getting homoclinic intersections.

We start the discussion by considering any partially hyperbolic set $\Lambda$
with a one-dimensional central bundle and whose strong bundles $E^{ss}$
and $E^{uu}$ are non trivial. One chooses a small cone field $\cC^c$ on a neighborhood of $\Lambda$
around the bundle $E^c$. This can be done in the following way:
one first extends continuously the subbundle $E^c$ in a neighborhood of $\Lambda$, one chooses
a small constant $\chi>0$ and one then defines at any point $x\in M$ close to $\Lambda$:
$$\cC^c_x=\{v\in T_xM, \; \|v^c\|> (1-\chi). \|v\|\},$$
where $v^c$ is the orthogonal projection of $v$ on $E^c_x$.\\
The open sets $\cC^c_x\setminus\{0\}$
have two components and for $r_0>0$ small enough,
the open balls $V$ in $M$ of radius $r_0$ are simply connected.
Hence, one can endow any of these balls with an orientation of $\cC^c_V$
(using for example a linear form which does not vanish on non-zero vectors of $\cC^c_V$.)\\
We also fix a small constant $L>0$.

Let $V$ be a ball of radius $r_0$ endowed with an orientation of the central cone field.
For any points $p,q\in K\cap V$, one says that
\emph{$p$ is below (resp. above) $q$} if there exists a positively
(resp. negatively) oriented $C^1$ curve tangent to $\cC^c$ of length smaller than $L$
that joins a point in the local strong unstable
manifold of $p$ to a point contained in the local strong stable manifold of $q$.
(One says that $p$ is \emph{strictly} below $q$ if the curve has non-zero length).

We mention two properties that justify these definitions (but that
will not be used in the following).
For points close enough, the position of $p$ with respect to $q$ is well defined:

\begin{lemma}
If $L$ is small enough, $p$ can not be both
strictly below and strictly above $q$.
\end{lemma}
\begin{proof}
It is similar to the proof of lemma 6.1 in \cite{bonatti-gan-wen}:
in small charts, the bundles $E^{ss}$, $E^{uu}$ and
the cone field $\cC^c$ are almost constant. If there exist two positively oriented
$C^1$ curves $[A,B]$ and $[C,D]$ tangent to $\cC^c$ with $A,D\in W^{ss}_{loc}(P)$
and $B,C\in W^{uu}_{loc}(Q)$, the sum $\overrightarrow{AB}+\overrightarrow{CD}$
belongs to a central cone field (since both vectors belong to
a central cone field and are positively oriented)
and the sum $\overrightarrow{AD}+\overrightarrow{CB}$ belongs to a sharp cone
around the space $E^{ss}\oplus E^{uu}$. This is a contradiction since
these sums coincide.
\end{proof}

On the other hand, any two close points $p,q\in K\cap V$ can always be compared:
\begin{lemma}
If $p,q$ are close enough, then $p$ is either below or above $q$.
\end{lemma}
\begin{proof}
It is similar to the proof of lemma~\ref{l.intermediate}:
one can cover a neighborhood of the ball $V$ with a smooth one-dimensional
foliation tangent to the central cone field $\cC^c$ such that any point in $V$ projects
by the holonomy along this foliation on a transverse one-codimensional disk $D$.
For any points $p,q\in V$, the local strong unstable manifold of $P$ and
the local strong stable manifold of $q$ intersect a leaf of the foliation, providing us with
a $C^1$-curve tangent to $C^c$ that joins the two local strong manifolds.
\end{proof}

\begin{figure}[ht]
\begin{center}
\input{twist.pstex_t}
\end{center}
\caption{Two points $p$, $q$ in a twisted position.\label{f.twist}}
\end{figure}
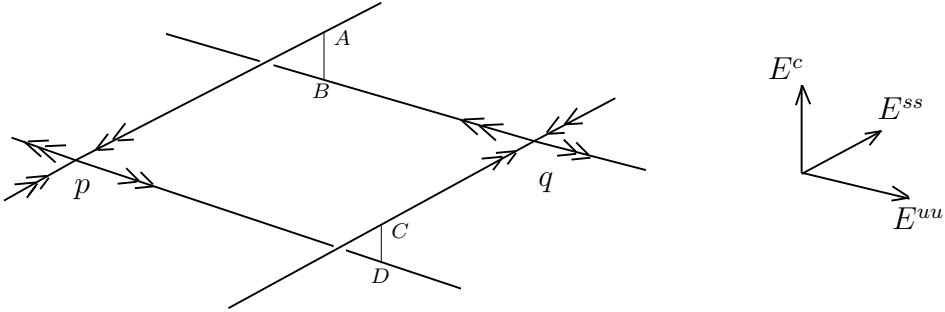
The fact that $p$ is below $q$ does not imply that $q$ is above $p$
(even if $p,q$ are close from each other).
This motivates the following definition (see figure~\ref{f.twist}):
\begin{definition}\label{d.twist}
Two points $p,q\in K\cap V$ are in a \emph{twisted position} if
\begin{itemize}
\item either $p$ is below $q$ and $q$ is below $p$,
\item or $p$ is above $q$ and $q$ is above $p$.
\end{itemize}
\end{definition}
Note that this definition does not depend on the choice of an orientation of
the central cone field on $V$.
One deduces that it may be extended for any pair of points $p,q\in \Lambda$ that are
close enough, without specifying a ball $V$ that contains both of them.

\begin{definition}\label{d.twist-returns}
We say that a sequence of periodic orbits $(\cO_n)_{n\in \NN}$ in $\Lambda$
has \emph{twisted returns} if there exists $\varepsilon>0$ small,
such that any points $p,q$ contained in a same orbit $\cO_n$ and at a distance smaller than $\varepsilon$
are in a twisted position.
\end{definition}
The definition remains unchanged if one replaces the cone field $\cC^c$ by a central cone field
defined on a smaller neighborhood of $\Lambda$ and which has sharper cones around $E^c$.
Since $\varepsilon>0$ can be chosen arbitrarily small, the definition
does not depend on the choice of a Riemannian metric nor of a central cone field.

Here is the perturbation result proved by Bonatti, Gan and Wen (we give a slightly
more precise statement than theorem 9 of \cite{bonatti-gan-wen}).

\begin{theorem}\label{t.twist}
Let $f$ be a diffeomorphism and
$\Lambda$ a partially hyperbolic set with a one-di\-men\-sio\-nal central bundle and
with non-trivial extremal bundles $E^{ss}$ and $E^{uu}$. Let $(\cO_n)$ be a sequence of periodic orbits contained in $\Lambda$ whose
periods go to infinity.
 
 If $(\cO_n)$ has twisted returns, then, for any neighborhood $\cU$ of $f$
in $\diff^1(M)$ and any neighborhood $U$ of $\Lambda$, there exist:
\begin{itemize}
\item a periodic orbit $\cO_{n_0}$ in the family $(\cO_n)$,
\item a diffeomorphism $g\in \cU$ that coincides with $f$ outside $U$ and on a
neighborhood of $\cO_{n_0}$,
\end{itemize}
such that the relative homoclinic class of $\cO_{n_0}$ in $U$ is non-trivial.
\end{theorem}

We explain the idea of the proof (see figure~\ref{f.twist}):
for a periodic orbit $\cO_{n_0}$ with a large period,
we consider in a chart two points $p$ and $q=f^k(p)$ of $\cO_{n_0}$
that minimize the distance $\dd(p,q)$ inside the finite set $\cO_{n_0}$.
They are in a twisted position and there
exist two small curves $[A,B]$ and $[C,D]$ tangent to the central cone
field such that $A$ and $C$ belong to the local strong stable manifold
of $p$ and $q$ respectively,
$D$ and $B$ belong to the local strong unstable manifold
of $p$ and $q$ respectively.
The bundles $E^{ss}$ and $E^{uu}$ in the chart are almost constant,
so that the twisted position implies that the lengths of 
$[A,B]$ and $[C,D]$ are very small in comparison to the distance
between $p$ and $q$. On the other hand,
since $\dd(p,q)$ is minimal, the four distances
$\dd(A, p)$, $\dd(D, p)$, $\dd(B,q)$ and
$\dd(C,q)$ are comparable to $\dd(p,q)$.
This shows that for any constant $\rho\geq 1$, one can assume that
the ball $V_\rho$ centered at $A$ and of radius
$\rho \dd(A,B)$ does not intersect $p$ nor $q$.
It does not intersect either any point $f^j(p)$ of the orbit $\cO_{n_0}$
since some forward iterate $f^{j+r}(p)$ would be close to
$f^r(A)$ and hence at a distance to $f^{r}(q)$ smaller
the distance $\dd(p,q)$.
After choosing a suitable constant $\rho>0$,
one can thus apply Hayashi's connecting lemma
in $V_\rho$ in order to join the backward orbit of $B$ to the forward
orbit of $A$ and create an intersection between the local strong
stable and unstable manifold of $p$ and $q$.

\begin{remark}
The proof produces an intersection between the strong stable and the strong unstable
manifolds of $\cO_{n_0}$. We already noticed (addendum~\ref{a.cycle})
that if the central Lyapunov exponent
of this periodic orbits is weak, then one can create a heterodimensional
cycle by a small $C^1$-perturbation of the dynamics.
\end{remark}

\begin{proof}[Proof of proposition~\ref{p.diff-strip} in the orientable case]$\quad$\\
For $\eta>0$ small, we denote by $D^{cs}(P_{\hat x})$ the ball of radius $\eta$
centered at $P_{\hat x}$ in the plaque $\cD_x^{cs}$ of $x=\pi(\hat x)$.
The local strong stable manifold $W^{ss}_{loc}(P_{\hat x})$ cuts this ball in two open disjoint components: $D^{cs}_-(P_{\hat x})$
and $D^{cs}_+(P_{\hat x})$ such that the latter is both contained in
$\hat \cD^{cs}_{\hat x}$ and disjoint from $\gamma_{\hat x}$.
We choose an orientation of $E^c$ on $K$ so that the parameterized curve $\pi(\{\hat x\}\times [0,+\infty))$
for each $\hat x\in \hat K$ is positively oriented at $x=\pi(\hat x)$.
By addendum \ref{a.model}.(\ref{i.model2}),
one can assume that the orientation of $E^c$ on $K$ extends
on $\Lambda=\pi(\hat\La)$ where $\hat \La$ is the maximal invariant set of $\hat f$ in $S$.
For a periodic point $\hat x$, the parameterized curve $\pi(\{\hat x\}\times [0,+\infty))$
first intersects $D^{cs}_-(P_{\hat x})$ and then $D^{cs}_+(P_{\hat x})$.

One can reduce to the case the bundle $E^{ss}$ is non-trivial:
if one assumes that $E^{ss}$ is trivial, the set $D^{cs}_{+}(P_{\hat x})$
is a curve of length $\eta>0$.
Let us denote by $B_{\hat x}$ the half ball
of radius $\eta$ centered at $P_{\hat x}$ in $M$, bounded by the (one-codimensional) local strong
unstable manifold $W^{uu}_{loc}(P_{\hat x})$ and intersecting the curve $D^{cs}_{+}(P_{\hat x})$;
its volume is uniformly bounded from below.
As in the proof of corollary~\ref{c.strong-unstable}, the period of $P_{\hat x}$
can be assumed arbitrarily large and there exist two distinct iterates $P_{\hat f^k(\hat x)}$
and $P_{\hat f^\ell (\hat x)}$ of $P_{\hat x}$ that are arbitrarily close. Their local strong unstable manifolds
are disjoint and have uniform size. Hence, (up to exchanging $k$ and $\ell$), the point $P_{\hat f^\ell (\hat x)}$
belongs to $B_{\hat f^k(\hat x)}$ and the local strong unstable manifold of $P_{\hat f^\ell (\hat x)}$
intersects transversally the curve $D^{cs}_{+}(P_{\hat f^k(\hat x)})$ at a point $z$.
By lemma~\ref{l.periodic}, $D^{cs}_{+}(P_{\hat f^k(\hat x)})$ is contained in the local stable
manifold of $P_{\hat f^k(\hat x)}$. The orbit of $z$ stays in an arbitrarily
small neighborhood $U$ of $\pi(S)$, hence, the relative homoclinic class of
$P_{\hat x}$ in $U$ is non-trivial, proving proposition~\ref{p.diff-strip} in this case.

By corollary~\ref{c.strong-unstable}, the bundle $E^{uu}$ also is non-trivial.
Hence, we can apply the beginning of section~\ref{sss.orientable} to the partially hyperbolic set $\Lambda$.
Let us fix a small constant $\varepsilon>0$.
If some point $P_{\hat x}$ has a return
$P_{\hat f^k( \hat x)}\neq P_{\hat x}$ at distance less than $\varepsilon$,
the local (strong) unstable manifold of $P_{\hat f^k(\hat x)}$ intersects transversally the ball
$D^{cs}(P_{\hat x})$. The local strong unstable manifold of $P_{\hat f^k(\hat x)}$
may thus be joined to the strong stable manifold of $W^{ss}_{loc}(P_{\hat x})$ by a small curve
tangent to a central cone $\cC^c$ and contained in one of the two half balls
$D^{cs}_-(P_{\hat x})$ or $D^{cs}_+(P_{\hat x})$.
By our choice of the global orientation on $\Lambda$, either the local unstable manifold
of $P_{\hat f^k(\hat x)}$ intersects $D^{cs}_-(P_{\hat x})$ and $P_{\hat f^k(\hat x)}$ is below
$P_{\hat x}$, or it intersects $D^{cs}_+(P_{\hat x})$
and $P_{\hat f^k(\hat x)}$ is above $P_{\hat x}$.
Since $D^{cs}_+(P_{\hat x})$ is contained in the local stable manifold of $P_{\hat x}$ (by lemma \ref{l.periodic}), the second case implies that the invariant manifolds of $P_{\hat x}$
intersect transversally at a point $z$. The orbit of $z$ stays in a small neighborhood of $\pi(S)$ and
proposition~\ref{p.diff-strip} again holds.

We are thus led to consider the situation where for any point $P_{\hat x}$
and any returns $P_{\hat f^k( \hat x)}$ and $P_{\hat f^\ell(\hat x)}$
that are $\varepsilon$-close,
$P_{\hat f^k(\hat x)}$ is below $P_{\hat f^\ell(\hat x)}$.
By symmetry, $P_{\hat f^\ell(\hat x)}$ is below $P_{\hat f^k(\hat x)}$,
so these points are in a twisted position and the sequence of orbits of the points
$P_{\hat x}$ has twisted returns. We already noticed that their periods go to infinity.
So, one can apply theorem~\ref{t.twist} and conclude that some periodic orbit of $\pi(S)$
has a non-trivial relative homoclinic class in $U$ after a $C^1$-small perturbation of $f$.
\end{proof}

%%%%%%%%%%%%%%%%%%%%%%%%%%%%%%%%%%%%%%%%%%%%%%%%%%%%%%%%%%%%%%%%%%%%%%%%%%%%%%%%
%%%%%%%%%%%%%%%%%%%%%%%%%%%%%%%%%%%%%%%%%%%%%%%%%%%%%%%%%%%%%%%%%%%%%%%%%%%%%%%%
\section{Existence of non-trivial homoclinic classes}
%%%%%%%%%%%%%%%%%%%%%%%%%%%%%%%%%%%%%%%%%%%%%%%%%%%%%%%%%%%%%%%%%%%%%%%%%%%%%%%%
\subsection{Some simple cases}\label{ss.simple}
Let $f$ be a $C^r$ diffeomorphism, $r\geq 1$.
We first recall how to get a transverse homoclinic intersection by a $C^r$-perturbation of $f$
in some particular cases. This shows that the closure of $\cM\cS\cup\cI$ in $\diff^r(M)$ contains
$\cH\text{yp}\cup\cT\text{ang}\cup\cC\text{yc}$.

\begin{description}
\item[1) The hyperbolic diffeomorphisms.] Let $f\in \cH\text{yp}$ be a hyperbolic diffeomorphism.
Since the Kupka-Smale diffeomorphisms are dense in $\diff^r(M)$ and since $\cH\text{yp}$ is open,
one can replace $f$ by a Kupka-Smale diffeomorphism $g$ arbitrarily $C^r$-close to $f$ and whose
chain-recurrent set $\cR(g)$ is hyperbolic. If $\cR(g)$ is finite, then $g$ is a Morse-Smale diffeomorphism.
In the other case $\cR(g)$ has a non-isolated point $z$. By the shadowing lemma, the periodic orbits
are dense in $\cR(g)$ and $z$ is the limit of a sequence of distinct hyperbolic periodic points.
The hyperbolicity of $\cR(g)$ then implies that all the periodic orbits that intersect a small neighborhood
of $z$ are homoclinically related. One deduces that $g$ has a transverse homoclinic intersection.
Thus $\left(\cM\cS\cup \cI\right)\cap \cH\text{yp}$ is dense in $\cH\text{yp}$. In particular,
$$\cH\text{yp}\subset \closure(\cM\cS\cup\cI).$$
\item[2) The homoclinic tangencies.] Let $f\in \cT\text{ang}$ be a diffeomorphism having a hyperbolic periodic point $p$
whose invariant manifolds have a non-transverse intersection point $z$.
By an arbitrarily $C^r$-small perturbation of $f$ in a neighborhood of $f^{-1}(z)$ one can modify the stable
manifold of $p$ at $z$ without changing the unstable manifold and create a transverse homoclinic intersection point. Hence,
$$\cT\text{ang}\subset \closure(\cI).$$
\item[3) The heterodimensional cycles.]
Let $f\in \cC\text{yc}$ be a diffeomorphism having two hyperbolic periodic orbits $\cO_1$ and $\cO_2$
such that the dimension of the stable spaces of $\cO_1$ is strictly less than the dimension of the stable
spaces of $\cO_2$ and such that moreover, the unstable manifold of $\cO_1$ (resp. of $\cO_2$)
intersects the stable manifold of $\cO_2$ (resp. of $\cO_1$) at a point $x$ (resp. $y$).
Since $\dim(W^u(\cO_1))+\dim(W^s(\cO_2))$ is larger than the dimension of $M$, one can,
up to a $C^r$-small perturbation of $f$, assume that the manifolds $W^u(\cO_1)$ and $W^s(\cO_2)$
are in general position at $x$. In particular, there exists a disk $D$ contained in $W^u(\cO_1)$
which intersects $W^s(\cO_2)$ transversally at $x$. By the $\lambda$-lemma, there exists a sequence of disks
$(D_n)_{n\geq 0}$ contained in $D$ such that the iterates $\left(f^n(D_n)\right)$ converge
towards a disk $D'\subset W^u(\cO_2)$ that contains $y$.
Consequently, there exist a small neighborhood $U$ of $y$ and a sequence $(y_n)$ converging towards $y$
such that the backward orbit of each $f^{-1}(y_n)$ is disjoint of $U$ and contained in $W^u(\cO_1)$.
By a $C^r$-small perturbation of $f$ in $f^{-1}(U)$, one thus can create a transverse homoclinic intersection point
for $\cO_1$. Hence,
$$\cC\text{yc}\subset \closure(\cI).$$
\end{description}

%%%%%%%%%%%%%%%%%%%%%%%%%%%%%%%%%%%%%%%%%%%%%%%%%%%%%%%%%%%%%%%%%%%%%%%%%%%%%%%%
\subsection{Existence of partially hyperbolic sets:
proof of theorem~A}
In this section we study the \emph{minimally non-hyperbolic sets},
i.e. the invariant compact sets that are non-hyperbolic and whose proper invariant compact subsets
are all hyperbolic. Zorn's lemma implies that the chain-recurrent set of
a non-hyperbolic diffeomorphism always contains minimally
non-hyperbolic sets. More precisely, we have:
\begin{proposition}\label{p.minimal}
Let $f$ be a diffeomorphism in $\cG_{rec}\setminus \cH\text{yp}$. Then, either $f$ belongs to $\cI$ or there exists a non-hyperbolic invariant set $K_0$ which carries a minimal dynamics.
\end{proposition}
\begin{proof}
Note that if all the chain-recurrence classes of $f$ are hyperbolic, the
chain-recurrent set $\cR(f)$ is hyperbolic, which is not the case since $f$ does not
belong to $\cH\text{yp}$. Hence, there exists a chain-transitive set $K_0$ of $f$ which is
non-hyperbolic. By Zorn's lemma, one can assume that $K_0$ is minimal for the inclusion.

If $K_0$ is not a minimal set, it contains a hyperbolic set $K\varsubsetneq K_0$.
By the shadowing lemma, the chain-recurrence class containing $K$ and $K_0$
contains a periodic orbit and (using that $f$ belongs to $\cG_{rec}$)
by theorem \ref{t.rec}, is a (non-trivial) homoclinic class.
This shows that $f$ has a transverse homoclinic orbit.
\end{proof}

In his work~\cite{wen} on the strong Palis conjecture,
Wen gets the following result (see~\cite[theorem B]{wen})
for transitive minimally non-hyperbolic sets.
\begin{theorem}\label{t.minimal}
There exists a dense G$_\delta$ subset $\cG_{mini}\subset \diff^1(M)\setminus \closure(\cT\text{ang}\cup\cC\text{yc})$
such that for any diffeomorphism $f\in \cG_{mini}$,
any invariant compact set $K_0$ which is transitive and minimally non-hyperbolic is
also partially hyperbolic. Moreover,
\begin{itemize}
\item either the central bundle $E^c$ of the partially hyperbolic splitting on $K_{0}$ is one-dimensional,
\item or the central bundle is two-dimensional and decomposes
into a dominated splitting $E^c=E^c_1\oplus E^c_2$ as the sum of two one-dimensional bundles.
\end{itemize}
\end{theorem}
The proof uses here S. Liao's selecting lemma (proposition 3.7 in~\cite{wen}, see also~\cite{liao})
to analyze a dominated splitting $T_{K_0}M=E\oplus F$ over a transitive minimally non-hyperbolic
set $K_{0}$. 
If neither $E$ nor $F$ is uniform, two cases can occur:
\begin{itemize}
\item Either Liao's selecting lemma can be used; in this case, $K_0$ intersects
a homoclinic class $H(P)$ (see lemma 3.8 in~\cite{wen}).
Moreover the dimension of the stable space of $P$ is equal to the dimension of $E$.
\item Or Liao's selecting lemma can not be applied and $K_0$ contains a hyperbolic set
$A$ (see lemma 3.12 in~\cite{wen}). Moreover, if $f$ does not belong to $\closure(\cT\text{ang})$,
one can choose $A$ such that the dimension of its stable bundle is less or equal to the dimension of $E$.
\end{itemize}
With these arguments Wen proves as a step towards theorem~\ref{t.minimal}
the following proposition (this is a restatement of proposition 3.13 in~\cite{wen},
but the proof is the same):
\begin{proposition}\label{p.deux-central}
Let $f$ be a diffeomorphism
and $K_0$ be an invariant compact set which is transitive and minimally non-hyperbolic.
Let us assume that the tangent bundle $T_{K_0}M$ over $K_0$ decomposes into a
(non-trivial) dominated splitting $E\oplus F$.

If the bundle $E$ is not uniformly contracted and the bundle $F$ is not uniformly expanded, then,
$K_{0}$ intersects a non-trivial homoclinic class $H(P)$. Moreover, if $f$ belongs to $\diff^1\setminus \closure(\cT\text{ang})$,
the dimension of the stable set of $P$ can be chosen less or equal to the dimension of the bundle $E$.
\end{proposition}

Note that if $f$ also belongs to $\cG_{rec}$ in the previous proposition,
the set $K_{0}$ is contained in $H(P)$.
One thus deduces from theorem \ref{t.minimal} and proposition \ref{p.deux-central}:
\begin{corollary}\label{c.minimal}
Let $f$ be a diffeomorphism in $\cG_{mini}\cap\cG_{rec}$ having a transitive
minimally non-hyperbolic set $K_{0}$. Then,
\begin{itemize}
\item either $K_0$ intersects a non-trivial homoclinic class and $f$ belongs to $\cI$,
\item or $K_{0}$ is partially hyperbolic with a one-dimensional central bundle.
\end{itemize}
\end{corollary}
\begin{proof}
By theorem \ref{t.minimal}, the set $K_{0}$ is partially hyperbolic.
If the central bundle is not one-dimensional, there exists a dominated splitting
$T_{K_{0}}=\left(E^{ss}\oplus E^{c}_{1}\right)\oplus \left(E^{c}_{2}\oplus E^{uu}\right)$
into two bundles that are not uniformly contracted or expanded.
One can thus apply proposition \ref{p.deux-central} to this decomposition and
$f$ has a non-trivial homoclinic class. Hence, $f$ belongs to $\cI$.
\end{proof}

The theorem~A is now
a direct consequence of these results.
\begin{proof}[Proof of theorem~A]
By theorem~\ref{t.minimal}, the set $\cG_{mini}\cup\cT\text{ang}\cup\cC\text{yc}$ is dense in $\diff^1(M)$;
so by section~\ref{ss.simple}, the set
$\cG_{mini}\cup\cM\cS\cup\cI$ has also this property.
Hence, by Baire theorem, the set
$$\cG_A=\left(\cG_{mini}\cup\cM\cS\cup\cI\right)\cap \cG_{rec}\cap \cG_{KS}$$
is a dense G$_\delta$ subset of $\diff^1(M)$.

Any diffeomorphism $f\in \cG_A\setminus \left(\cM\cS\cup\cI\right)$ belongs
to $\cG_{rec}\setminus \cH\text{yp}$ (by section~\ref{ss.simple})
and by proposition~\ref{p.minimal} has a non-hyperbolic invariant compact set $K_0$ which carries a minimal dynamics. In particular, $K_0$ is a transitive
minimally non-hyperbolic set. Since $f$ is Kupka-Smale,
$K_0$ is not a periodic orbit.
Since $f$ belongs to $\left(\cG_{mini}\cap \cG_{rec}\right)\setminus \cI$,
one can apply corollary~\ref{c.minimal}: the set $K_0$ is partially hyperbolic with
a one-dimensional central bundle, as required.
\end{proof}

%%%%%%%%%%%%%%%%%%%%%%%%%%%%%%%%%%%%%%%%%%%%%%%%%%%%%%%%%%%%%%%%%%%%%%%%%%%%%
\subsection{Approximation of partially hyperbolic sets by homoclinic classes:
proof of theorem~B}
Theorem~B will be
a consequence of the following perturbation result:

\begin{proposition}\label{p.partial}
Let $f$ be a Kupka-Smale diffeomorphism in the dense G$_\delta$ subset $\cG_{rec}\cap\cG_{shadow}$
of $\diff^1(M)$ and $K_0$ be a chain-transitive set of $f$
which is partially hyperbolic with a one-dimensional central bundle.

If $K_0$ is not a periodic orbit, then, for any neighborhood $\cU$ of $f$ in $\diff^1(M)$, any neighborhood $U$ of $K_0$ and any $\varepsilon>0$, there exist:
\begin{itemize}
\item a periodic orbit $\cO$ that is $\varepsilon$-close to $K_0$ for the Hausdorff distance,
\item a diffeomorphism $g\in \cU$ that coincides with $f$ outside $U$ and on a neighborhood
of $\cO$, such that the relative homoclinic class of $\cO$ in $U$ is non-trivial for $g$.
\end{itemize}
\end{proposition}
\begin{proof}
If one assumes that $K_0$ contains a periodic orbit, then
(since $f$ belongs to $\cG_{rec}$) by theorem~\ref{t.rec} the set $K_0$ is contained in a
non-trivial relative homoclinic class $H(P,V)$ where $V\subset U$ is an arbitrarily small neighborhood of $K_{0}$. Hence, the conclusion of proposition~\ref{p.partial} is satisfied by the diffeomorphism $g=f$
and any periodic orbit $\cO$ related in $V$ to the orbit of $P$ and close to $H(P,V)$
for the Hausdorff distance.

We now assume that $K_0$ contains no periodic point.
Since $f$ belongs to $\cG_{shadow}$, by theorem~\ref{t.shadow} the set $K_{0}$ is the limit
for the Hausdorff distance of a sequence of periodic orbits $\left(\cO_n\right)$.
Since $K_0$ is not a single periodic orbit, the period of $\cO_n$ goes to infinity with $n$.
One introduces the set $K$ as the union of $K_0$ with the periodic orbits $\cO_n$ for $n\geq n_0$.
If $n_0$ is large enough, the set $K$ is arbitrarily close to $K_0$ for the Hausdorff distance.
In particular, it is a partially hyperbolic set with a one-dimensional central bundle.
If $f$ preserves an orientation on the central bundle $E^c$ over $K_0$, then
$f$ also preserves an orientation on the central bundle over $K$, by proposition~\ref{p.orientation}.

Now, we apply proposition~\ref{p.existence} and introduce a central model $(\hat K,\hat f)$
for $(K,f)$ with a projection $\pi$ which satisfies the properties stated in addendum~\ref{a.model}.
By considering the points in $\hat K$ that projects by $\pi$ on $K_0$,
one also obtains a central model $(\hat K_0,\hat f)$ for $(K_0,f)$.
Since $K_0$ is chain-transitive, the base of $(\hat K_0,\hat f)$ is chain-transitive:
\begin{itemize}
\item Either $f$ preserves an orientation on the central bundle $E^c$ on $K$. By addendum~\ref{a.model}.(\ref{i.model3}), the map
$\pi\colon \hat K\times \{0\}\to K$ is a bijection and the dynamics
of $\hat f$ on $\hat K_0\times \{0\}$ is conjugate to the dynamics of $f$ on $K_0$.
\item Or, $f$ does not preserve any orientation of $E^c$ on $K_0$ (nor on $K$). By addendum~\ref{a.model}.(\ref{i.model4}),
$\pi\colon \hat K_0\to K_0$ is the bundle of orientations of $E^c$ over $K_0$ and lemma~\ref{l.orientation}
implies that $(\hat K_0,\hat f)$ is chain-transitive.
\end{itemize}

If the central model $(\hat K_0,\hat f)$ has a chain-recurrent central segment,
one can apply proposition~\ref{p.segment} (using that $f$ belongs to $\cG_{rec}\cap \cG_{shadow}$):
in this case, $K_0$ is contained in a (non-trivial) relative homoclinic class $H(P,U)$
in $U$ and the conclusion of proposition~\ref{p.partial} is again satisfied.

We thus consider the case the central model $(\hat K_0,\hat f)$ has no chain-recurrent central segment:
by proposition~\ref{p.strip}, it has a trapping strip $S_0$
that is an arbitrarily small neighborhood of $\hat K_0\times \{0\}$.
By continuity of $\hat f$, if $\hat K$ is close enough to $\hat K_0$
(which can be supposed by increasing $n_0$),
this implies that $(\hat K,\hat f)$ also has a trapping strip $S$.
Note that any periodic orbit $\hat \cO$ of $\hat f$ should be contained in
$(\hat K\setminus \hat K_0)\times [0,+\infty)$ since $K_0$ has no periodic point.
As $K$ is the limit of the periodic orbits $\cO_n$, by assuming $n_0$ large enough
and the strip $S_0$ small, this shows that the projection of any periodic orbit
$\hat \cO$ of $\hat f$ contained in $S$ is arbitrarily close to $K_0$
for the Hausdorff distance.\\
Since $f$ is Kupka-Smale, one can apply proposition~\ref{p.diff-strip}: there exist
\begin{itemize}
\item a periodic orbit $\cO$ of $f$ which is the projection by $\pi$ of
a periodic orbit $\hat \cO\subset S$ of $\hat f$ (hence $\cO$ is
$\varepsilon$-close to $K_0$ for the Hausdorff distance),
\item a diffeomorphism $g\in \cU$ that coincides with $f$ outside $U$ and on a neighborhood
of $\cO$ such that the relative homoclinic class of $\cO$ in $U$ is non-trivial.
\end{itemize}
This ends the proof of the proposition.
\end{proof}

One concludes the proof of the theorem by a Baire argument.
\begin{proof}[Proof of theorem B]
Let us consider a dense countable subset $X$ of $M$. For any finite set $\{x_1,\dots,x_r\}$ in $X$,
and any $n\geq 1$, we define the set $\cI(x_1,\dots,x_r,n)$ of diffeomorphisms in $\diff^1(M)$
having a non-trivial relative homoclinic class whose Hausdorff distance
from the set $\{x_1,\dots,x_r\}$ is strictly less than $\frac{1}{n}$. This set is open.
We denote by $\cU(x_1,\dots,x_r,n)$ the set $\diff^1(M)\setminus \boundary(\cI(x_1,\dots,x_r,n))$
which is open and dense in $\diff^1(M)$.

By Baire theorem, the set
$$\cG_B=\cG_{KS}\cap \cG_{rec}\cap \cG_{shadow}\cap \bigcap_{(x_1,\dots,x_r,n)}\cU(x_1,\dots,x_r,n)$$
of Kupka-Smale diffeomorphism that satisfy both theorems~\ref{t.rec} and~\ref{t.shadow} and that
are contained in each $\cU(x_1,\dots,x_r,n)$ is a dense G$_\delta$ subset of $\diff^1(M)$.

For a diffeomorphism $f\in \cG_B$, we consider a chain-transitive partially hyperbolic set $K_0$
which has a one-dimensional central bundle and which is not a periodic orbit. We fix $n\geq 1$.
There exists a set $\{x_1,\dots,x_r\}$ whose Hausdorff distance to $K_{0}$
is strictly less than $\frac{1}{2n}$.
Since $f$ belongs to $\cG_{KS}\cap\cG_{rec}\cap \cG_{shadow}$, proposition~\ref{p.partial}
may be applied: there exists $C^1$-small perturbations of $f$
having a non-trivial relative homoclinic class at distance less than $\frac{1}{2n}$ from
$K_0$, hence at distance strictly less than $\frac 1 n$ from $\{x_{1},\dots,x_{r}\}$.
One deduces that $f$ is accumulated by diffeomorphisms in $\cI(x_1,\dots,x_r,n)$.
Since $f$ belongs to $\cU(x_1,\dots,x_r,n)$, it also belongs to $\cI(x_1,\dots,x_r,n)$.
Thus, $f$ has a non-trivial relative homoclinic class $\frac{1}{n}$-close to $\{x_{1},\dots,x_{r}\}$ and $\frac{3}{2n}$-close to $K_{0}$ for
the Hausdorff distance. This shows that $K_0$ is the limit of non-trivial relative homoclinic classes, as required.
\end{proof}

\section*{Appendix: the case of conservative and tame diffeomorphisms}

We first prove the strong Palis conjecture in the case of tame diffeomorphisms.
\begin{theorem*}
Any diffeomorphism of a compact manifold can be $C^1$-approximated
by a diffeomorphism which is hyperbolic or has a heterodimensional cycle
or admits filtrations with an arbitrarily large number of levels.
\end{theorem*}
\begin{proof}
We first perturb the diffeomorphism in $\diff^1(M)$ so that
it satisfies~\cite[theorem C]{flavio} and the following generic properties:
any homoclinic class is a chain-recurrence class and
any isolated chain-recurrence class in the chain-recurrent set is a homoclinic class (see~\cite{BC}).

Now, either the new diffeomorphism has infinitely many chain-recurrence classes,
or it has finitely many homoclinic classes only.
In the first case, the diffeomorphism admits filtrations with an arbitrarily large number of levels. In the second case,
\cite[theorem C]{flavio} concludes that this diffeomorphism is either
hyperbolic or can be $C^1$-approximated by a diffeomorphism which exhibits
a heterodimensional cycle.
\end{proof}
\bigskip

In the conservative setting, the argument is the same, but
we have to check that all the perturbations involved can be made
conservative.
\begin{theorem*}
Any conservative diffeomorphism of a compact surface can be $C^1$-approximated
by a conservative diffeomorphism which is Anosov or has a homoclinic tangency.

Any conservative diffeomorphism of a manifold $M$ of dimension
$\dim(M)\geq 3$ can be $C^1$-approximated by a conservative diffeomorphism which is Anosov or has a heterodimensional cycle.
\end{theorem*}
\begin{proof}
Using~\cite{ABC}, one can find a $C^1$-small conservative perturbation $f$
so that $M$ is the homoclinic class of a hyperbolic periodic point $P$.
We denote by $0<i<\dim(M)$ the stable dimension of $P$.
We may also assume (by~\cite{robinson-periodic} and~\cite[proposition 7.4]{bdp}) that there exists a $C^1$-neighborhood $\cU$
of $f$ in the space of conservative diffeomorphisms such that
one of the following properties holds:
\begin{enumerate}
\item[a)] either $f$ has a periodic point with complex eigenvalues (if $M$
is a surface) or which is hyperbolic with stable dimension different from $i$
(if $\dim(M)\geq 3$),
\item[b)] or the periodic points of any $g\in \cU$are hyperbolic
and have stable dimension equal to $i$.
\end{enumerate}

On a compact surface, if a conservative diffeomorphism has a periodic point $P$
with a complex eigenvalue, then one can create a homoclinic
tangency by a $C^1$-small conservative perturbation: using generating functions
one can first perturb so that the diffeomorphism at the period
is $C^1$ conjugated to an irrational rotation in a neighborhood of the periodic point; then, one creates close to $P$ a periodic disc such that the return map coincides with
the identity; by a last perturbation we obtain a homoclinic tangency
inside this disc as a perturbation of the identity with compact support.

In higher dimension and for a $C^1$-generic (conservative) diffeomorphism $f$,
if the homoclinic class of a periodic point $P$
contains a hyperbolic periodic point with different stable dimension,
then $f$ can be $C^1$-approximated by a diffeomorphism having a heterodimensional
cycle: this is a consequence of Hayashi's connecting lemma.

We conclude in case a) that $f$ can be approximated by a diffeomorphism
which has a homoclinic tangency or a heterodimensional cycle depending if $\dim(M)=2$ or $\geq 3$.
\medskip

In case b), one has to reproduce the arguments of Ma\~n\'e~\cite{mane}
in the conservative setting. First:
\begin{itemize}
\item[--] there exists on $M$ a dominated splitting
$T_M=E\oplus F$ with $\dim(E)=i$,
\item[--] there exists $N\geq 1$ and constants $C>0$, $\lambda\in (0,1)$
such that for any periodic orbit point $P$ of period $\tau$ one has
\begin{equation}\label{e.contraction}
\prod_{k=1}^\tau \|Df^N_{|E}(f^k(P))\|\leq C. \lambda^\tau
\text{ and }
\prod_{k=1}^\tau \|Df^{-N}_{|F}(f^k(P))\|\leq C. \lambda^\tau.
\end{equation}
\end{itemize}
Otherwise by~\cite[lemma II.3]{mane}
there exists a small perturbation
of the differential $Df$ above one of the periodic orbits of $f$ which
has an eigenvalue of modulus one.
This perturbation can be taken conservative:
it is enough to compose the tangent map by linear maps which are the identity on the central
spaces and which are a homotheties on their orthogonal spaces.
With~\cite[proposition 7.4]{bdp} one then creates
a perturbation of $f$ with a non-hyperbolic periodic point,
contradicting the assumption b).

We then prove that the bundle $E$ is uniformly contracted
and more precizely that for any invariant probability measure $\mu$
one has
$$\int \log \|Df^N_{|E}\|d\mu<0.$$
If this fails, by the ergodic closing lemma~\cite{mane}
(for the conservative setting) there exists a sequence of periodic
points $(P_n)$ of $f$ with periods $(\tau_n)$
such that 
$$\frac 1 {\tau_n} \sum_{k=1}^{\tau_n} \log \|Df^N_{|E}(f^k(P_n))\|$$
converges to a non-negative number, contradicting~(\ref{e.contraction}).
The same argument shows that $F$ is uniformly expanded.
We have thus proved that $f$ is an Anosov diffeomorphism.
\end{proof}

%%%%%%%%%%%%%%%%%%%%%%%%%%%%%%%%%%%%%%%%%%%%%%%%%%%%%%%%%%%%%%%%%%%%%%%%%%%%%%

\vskip 1cm

\flushleft{\bf Sylvain Crovisier} \\
CNRS - Laboratoire Analyse, G\'eom\'etrie et Applications, UMR 7539,\\
Institut Galil\'ee, Universit\'e Paris 13, Avenue J.-B. Cl\'ement, 93430 Villetaneuse, France\\
\textit{E-mail:} \texttt{crovisie@math.univ-paris13.fr}
\end{document}

%% file: central.pstex_t
\begin{picture}(0,0)%
\includegraphics{central.pstex}%
\end{picture}%
\setlength{\unitlength}{1579sp}%
\begingroup\makeatletter\ifx\SetFigFont\undefined%
\gdef\SetFigFont#1#2#3#4#5{%
  \reset@font\fontsize{#1}{#2pt}%
  \fontfamily{#3}\fontseries{#4}\fontshape{#5}%
  \selectfont}%
\fi\endgroup%
\begin{picture}(6690,9017)(2686,-9281)
\put(8401,-9136){\makebox(0,0)[b]{\smash{{\SetFigFont{10}{12.0}{\rmdefault}{\mddefault}{\updefault}{$\gamma$}%
}}}}
\put(6526,-1036){\makebox(0,0)[b]{\smash{{\SetFigFont{10}{12.0}{\rmdefault}{\mddefault}{\updefault}{$P$}%
}}}}
\put(8506,-1411){\makebox(0,0)[b]{\smash{{\SetFigFont{10}{12.0}{\rmdefault}{\mddefault}{\updefault}{$B$}%
}}}}
\put(8506,-8311){\makebox(0,0)[b]{\smash{{\SetFigFont{10}{12.0}{\rmdefault}{\mddefault}{\updefault}{$A$}%
}}}}
\put(9256,-3136){\makebox(0,0)[b]{\smash{{\SetFigFont{10}{12.0}{\rmdefault}{\mddefault}{\updefault}{$\pi$}%
}}}}
\put(8536,-5596){\makebox(0,0)[b]{\smash{{\SetFigFont{10}{12.0}{\rmdefault}{\mddefault}{\updefault}{$z$}%
}}}}
\put(4456,-9031){\makebox(0,0)[b]{\smash{{\SetFigFont{10}{12.0}{\rmdefault}{\mddefault}{\updefault}{$L$}%
}}}}
\put(4261,-2701){\makebox(0,0)[b]{\smash{{\SetFigFont{10}{12.0}{\rmdefault}{\mddefault}{\updefault}{$q$}%
}}}}
\put(4246,-1396){\makebox(0,0)[b]{\smash{{\SetFigFont{10}{12.0}{\rmdefault}{\mddefault}{\updefault}{$p$}%
}}}}
\put(5986,-3721){\makebox(0,0)[b]{\smash{{\SetFigFont{10}{12.0}{\rmdefault}{\mddefault}{\updefault}{$Q$}%
}}}}
\put(2701,-5836){\makebox(0,0)[b]{\smash{{\SetFigFont{10}{12.0}{\rmdefault}{\mddefault}{\updefault}{$D$}%
}}}}
\end{picture}%

%% file: periodic.pstex_t
\begin{picture}(0,0)%
\includegraphics{periodic.pstex}%
\end{picture}%
\setlength{\unitlength}{1579sp}%
\begingroup\makeatletter\ifx\SetFigFont\undefined%
\gdef\SetFigFont#1#2#3#4#5{%
  \reset@font\fontsize{#1}{#2pt}%
  \fontfamily{#3}\fontseries{#4}\fontshape{#5}%
  \selectfont}%
\fi\endgroup%
\begin{picture}(8525,8481)(-691,-8179)
\put(5071,-4336){\makebox(0,0)[b]{\smash{{\SetFigFont{12}{14.4}{\rmdefault}{\mddefault}{\updefault}{$x$}%
}}}}
\put(7801,-961){\makebox(0,0)[lb]{\smash{{\SetFigFont{12}{14.4}{\rmdefault}{\mddefault}{\updefault}{$D_+^ {cs}(P_{\hat x})$}%
}}}}
\put(5101,-46){\makebox(0,0)[b]{\smash{{\SetFigFont{12}{14.4}{\rmdefault}{\mddefault}{\updefault}{$\sigma_{\hat x}$}%
}}}}
\put(7816,-5146){\makebox(0,0)[lb]{\smash{{\SetFigFont{12}{14.4}{\rmdefault}{\mddefault}{\updefault}{$\cD^{cs}_x$}%
}}}}
\put(5326,-2086){\makebox(0,0)[b]{\smash{{\SetFigFont{12}{14.4}{\rmdefault}{\mddefault}{\updefault}{$P_{\hat x}$}%
}}}}
\put(1816,-4081){\makebox(0,0)[rb]{\smash{{\SetFigFont{8}{9.6}{\rmdefault}{\mddefault}{\updefault}{$W^{ss}_{loc}(x)$}%
}}}}
\put(1846,-1111){\makebox(0,0)[rb]{\smash{{\SetFigFont{8}{9.6}{\rmdefault}{\mddefault}{\updefault}{$\widehat\cD^ {cs}_{\hat x}$}%
}}}}
\put(1786,-3076){\makebox(0,0)[rb]{\smash{{\SetFigFont{8}{9.6}{\rmdefault}{\mddefault}{\updefault}{$\gamma_{\hat x}$}%
}}}}
\end{picture}%

%% file: twist.pstex_t
\begin{picture}(0,0)%
\includegraphics{twist.pstex}%
\end{picture}%
\setlength{\unitlength}{1579sp}%
\begingroup\makeatletter\ifx\SetFigFont\undefined%
\gdef\SetFigFont#1#2#3#4#5{%
  \reset@font\fontsize{#1}{#2pt}%
  \fontfamily{#3}\fontseries{#4}\fontshape{#5}%
  \selectfont}%
\fi\endgroup%
\begin{picture}(15311,5064)(43,-6463)
\put(5986,-5901){\makebox(0,0)[b]{\smash{{\SetFigFont{8}{9.6}{\rmdefault}{\mddefault}{\updefault}{$D$}%
}}}}
\put(5056,-2991){\makebox(0,0)[b]{\smash{{\SetFigFont{8}{9.6}{\rmdefault}{\mddefault}{\updefault}{$B$}%
}}}}
\put(6301,-5181){\makebox(0,0)[b]{\smash{{\SetFigFont{8}{9.6}{\rmdefault}{\mddefault}{\updefault}{$C$}%
}}}}
\put(5386,-2151){\makebox(0,0)[b]{\smash{{\SetFigFont{8}{9.6}{\rmdefault}{\mddefault}{\updefault}{$A$}%
}}}}
\put(1306,-4476){\makebox(0,0)[b]{\smash{{\SetFigFont{12}{14.4}{\rmdefault}{\mddefault}{\updefault}{$p$}%
}}}}
\put(8581,-4371){\makebox(0,0)[b]{\smash{{\SetFigFont{12}{14.4}{\rmdefault}{\mddefault}{\updefault}{$q$}%
}}}}
\put(12346,-2691){\makebox(0,0)[b]{\smash{{\SetFigFont{12}{14.4}{\rmdefault}{\mddefault}{\updefault}{$E^c$}%
}}}}
\put(14431,-5046){\makebox(0,0)[b]{\smash{{\SetFigFont{12}{14.4}{\rmdefault}{\mddefault}{\updefault}{$E^{uu}$}%
}}}}
\put(14161,-3321){\makebox(0,0)[b]{\smash{{\SetFigFont{12}{14.4}{\rmdefault}{\mddefault}{\updefault}{$E^{ss}$}%
}}}}
\end{picture}%